\documentclass[msom,nonblindrev]{informs3}

\DoubleSpacedXI

\usepackage{amsmath}
\usepackage{amssymb}
\usepackage{algorithm}
\usepackage{algpseudocode}
\usepackage{caption}
\usepackage{subcaption}
\usepackage{graphicx}
\theoremstyle{THkey}\newtheorem{assumptiont}{XXXXX}

\usepackage{natbib}
 \bibpunct[, ]{(}{)}{,}{a}{}{,}%
 \def\bibfont{\small}%
\TheoremsNumberedThrough     

\EquationsNumberedThrough    

\renewcommand{\theARTICLETOP}{}

\begin{document}


\RUNAUTHOR{Lee et al.}

\RUNTITLE{Optimal Policy for Inventory Management with Periodic and Controlled Resets}

\TITLE{Optimal Policy for Inventory Management with Periodic and Controlled Resets}

\ARTICLEAUTHORS{%
\AUTHOR{Yoon Lee}
\AFF{Department of Industrial Engineering and Operations Research, University of California, Berkeley, CA 94720, \EMAIL{yllee@berkeley.edu}}
\AUTHOR{Yonatan Mintz}
\AFF{Department of Industrial and Systems Engineering, University of Wisconsin-Madison, WI 53706, \EMAIL{ymintz@wisc.edu}}
\AUTHOR{Anil Aswani, Zuo-Jun Max Shen}
\AFF{Department of Industrial Engineering and Operations Research, University of California, Berkeley, CA 94720, \EMAIL{\{aaswani,maxshen\}@berkeley.edu}}
\AUTHOR{Cong Yang}
\AFF{Sauder School of Business, University of British Columbia, Vancouver, BC V6T 1Z2, Canada, \EMAIL{cong.yang@sauder.ubc.ca}}
} 

\ABSTRACT{%
Inventory management problems with periodic and controllable resets occur in the context of managing water storage in the developing world and retailing limited-time availability products. In this paper, we consider a set of sequential decision problems in which the decision-maker must not only balance holding and shortage costs but discard all inventory before a fixed number of decision epochs, with the option for an early inventory reset. Finding optimal policies using dynamic programming for these problems is particularly challenging since the resulting value functions are non-convex. Moreover, this structure cannot be easily analyzed using existing extended definitions, such as $K$-convexity. Our key contribution is to present sufficient conditions that ensure the optimal policy has an easily interpretable structure that generalizes the well-known $(s, S)$ policy from the operations literature. Furthermore, we demonstrate that the optimal policy has a four-threshold structure under these rather mild conditions. We then conclude with computational experiments, thereby illustrating the policy structures that can be extracted in several inventory management scenarios.
%
%
}%


\KEYWORDS{supply chain management, inventory theory and control, dynamic programming, health care management, humanitarian operations}

\maketitle

%


\section{Introduction}

The study of inventory management with stochastic demand has been a key discipline in operations management since the inception of the field. In general, these problems consider a decision-maker who must choose either an order quantity or a replenishment level such that it minimizes their overall expected supply chain costs or that they are able to maintain a certain level of service with high probability. Several models have been proposed in the literature to address this setting, such as the newsvendor model for single period replenishment \citep{arrow1951}, the $(r,Q)$ model for continuous review policies \citep{hadley1963}, and the $(s,S)$ model for periodic inventory review models \citep{scarf1959}. These settings generally focus on finding an optimal inventory policy, that is, a decision rule that maps the current inventory level of the facility to an order quantity. To find this policy, classical models have mainly focused on balancing holding costs (the costs incurred from holding excess inventory) and shortage costs (the costs associated with not having enough units on hand to satisfy demand). However, when products are perishable or subject to strict health and safety regulations, the assumptions that underlie these policies are no longer met.

Many settings, such as food procurement \citep{van2006, blackburn2009, farahani2012}, medical supply chain management \citep{prastacos1984, pierskalla2005, shen2011}, and certain seasonally sensitive retail supply chains \citep{caro2007, nagurney2013}, require managing perishable inventory and therefore do not satisfy the classical settings of the aforementioned inventory management models. The common way this constraint is addressed is through the use of strategic discarding policies. Essentially, the decision-maker is able to form a ranking of which products are most valuable to keep for longer durations and to individually discard expired units based on costs and health constraints. While these policies are applicable to many different contexts, if the particular goods being managed are commingled or not easily separable, it is not feasible to strategically discard individual units, but instead, the entirety of the inventory must be discarded. For example, in the case of residential water storage, as new water is purchased, it mixes with older water that is more likely to have been tainted and cannot be readily separated from the previous water in the tank.

In this paper, we propose new models and techniques to address the setting of a perishable inventory with known expiration time, where strategic discarding is not possible. In particular, we propose to model these systems as inventory control problems with strict reset constraints. In this setting, at each epoch, the decision-maker is faced with two decisions: first, whether or not to discard the entire inventory on hand; and second, whether or not to order inventory up to a certain level. Since the inventory is set to expire at a predetermined time, the decision-maker must discard the entirety of the inventory after a fixed number of epochs have passed. The main challenge in this problem is to find an inventory policy for the decision-maker that balances safety constraints, in conjunction with holding and shortage costs, by making refresh and reorder decisions.


\subsection{Applications}

The setting of mixing perishable inventory with reset decisions is common in many real-world applications. Here, we present two settings that exemplify the assumptions of the model. First, we consider a non-profit example of managing water storage in a residential building in the developing world. Next, we describe a profit example in the case of a retail supply chain with switching product lines.

\subsubsection{Non-Profit: Water Storage Problem}

Over 300 million people around the world have an intermittent water supply \citep{kumpel2016}, making water availability a pervasive public health crisis. This problem is more severe in the developing world, where the current distribution infrastructure is not capable of providing cities with continuous water supplies.
%
%
Consequently, most households resort to communal and personal water storage containers to maintain a water supply throughout the day, and these water storage containers are filled once every few days during intermittent periods of water availability \citep{worldbank2010, who2017}. However, water stored in these containers is prone to contamination, and water stored long-term is likely to have higher concentrations of bacterial and viral pathogens \citep{gadgil1998, zerah1998, coelho2003, tokajian2003, lee2005, klingel2012}. 
The problem of managing local water storage can be posed in the framework of perishable inventory that we have described previously. Each day, the decision-maker must choose how much additional water to purchase to satisfy the demand for that day. Since the longer water sits in the tank, the higher its chance of getting contaminated, this can be modeled as a holding penalty to ensure that water does not get stored for too long. Likewise, if the amount of water purchased is lower than the demand realized for that day, expediting the shipment of additional water may be costly, which can be modeled as a shortage penalty. Finally, since the risk of contamination makes water perishable and contaminated water cannot be easily separated from clean water, the only option available to the decision-maker is a ``reset'' type action where the whole amount of water on hand is discarded.
%

\subsubsection{Profit: Retail Management Problem}

In many retail supply chain settings, product demand is highly dependent on seasonal tastes or fashion trends. This can affect various types of products, ranging from apparel and food to consumer electronics \citep{matsuo2007, caro2015}. In the specific case of clothing, this has previously been studied in the operations literature through the lens of fast fashion \citep{barnes2006, caro2010, bhardwaj2010, cachon2011}. The key to managing these inventories is that as new trends arrive, new product lines must be introduced to meet the most current customer needs. While production and inventory costs may remain the same for new product lines, as products fall out of fashion with consumers, the old product line must be completely discarded to make room for the new product line. This setting can again be modeled by our perishable inventory with reset control framework. Although, unlike water, apparel and similar products are easy to separate, due to how tastes shift in terms of seasonality and fashion, all products become obsolete simultaneously and not individually. This means that strategic discarding is not possible since the entirety of the stock needs to be replaced. Hence, the key challenge to designing an inventory policy in this setting is similar to that of water storage. Not only do inventory ordering, shortage, and holding costs need to be balanced, but also the decision-maker must consider at what point it is advantageous to take a ``reset'' action corresponding to introducing a new product line.

\subsection{Related Literature}

In this paper, we examine the problem of perishable inventory with reset control, a special case of multi-period inventory management with periodic review and stochastic demand. In the context of supply chain theory, a policy is a function that maps the current state of our inventory into an ordering decision \citep{snyder2019}. For periodic-review inventory models with stochastic demands, a base-stock policy (e.g., newsvendor model) or an $(s, S)$ policy is an example of inventory control policy. The fundamental idea underlying the $(s, S)$ policy is the following: in each time period, we observe the current inventory position---if the inventory position falls below $s$, then we place an order of sufficient size to bring the inventory position to $S$, where the quantity $s$ is known as the reorder point and $S$ as the order-up-to level \citep{arrow1951, scarf1959}. To this classical framework, we seek to add a reset control action to empty the entire inventory and generalize the $(s, S)$ policy for such problems. \cite{mintz2017} studied this approach numerically, and our contributions include the theoretical analysis of the optimal policy in the presence of periodic and controllable resets. In addition to these works, we draw on two main streams of literature for the basis of our modeling and analysis, namely inventory management and dynamic programming.

\subsubsection{Inventory and Supply Chain Management}

The problem we describe in this paper is closely related to the setting of managing perishable inventory \citep{pierskalla1972, nahmias1982, hsu2000, hsu2001, chu2005, karaesmen2011, coelho2014}, which has been explored primarily in the context of healthcare settings. In these perishable inventory models, other than the cost of production, costs arise either from not having enough inventory to meet demand (shortage costs) or from holding excess inventory (holding costs) that must be discarded due to the perishable nature of the good. In contrast to our setting, this stream of literature generally assumes that expired units of inventory can be disposed of individually, which allows salvage and disposal costs to be incorporated into the holding cost, and hence, policies of strategic inventory removal can be employed \citep{rosenfield1989, rosenfield1992}. In our setting, however, decision-makers cannot dispose of single units of their stock strategically. For example, in the water storage problem, fresh water may mix with still water, and they cannot be easily separated. Instead, the decision-maker must decide to either keep or dispose of their entire inventory in order to reduce the risk of drinking contaminated water. Similarly, in the retail management problem, we determine if an entire product line should be kept or discontinued. In the case where the product line is discontinued, all of its inventory must be discarded to free up space for the inventory of the new product line.

\subsubsection{Dynamic Programming and Optimal Control}

It is well known from the operations literature that periodic-review inventory management problems can be modeled as continuous-state dynamic programs \citep{scarf1959,bertsekas1995}. For most of these models, while it is possible to derive structural results for the functional form of the optimal policy, it is often difficult to find a closed-form solution for the parameters of the policy. This means that, in general, numerical algorithms must be used to find the relevant parameter values. The main two numerical methods developed to solve these dynamic programming problems are value and policy iteration \citep{bertsekas1995}, which require finding a fixed point in either the value function or the policy space, respectively. Essentially, this means algorithms are attempting to find fixed points in functional spaces that may be infinite dimensional. Even in finite high-dimensional state spaces, the exact calculation for these approaches is numerically difficult, an effect commonly referred to as the \emph{curse of dimensionality}. To address this challenge, a number of approximate dynamic programming approaches \citep{bertsekas1995,de2003linear,powell2007approximate,ryzhov2012knowledge,kariotoglou2013approximate,haskell2016empirical} have been developed. The key to these approaches is to perform a value or policy iteration with an inexact, but tractable, representation of the value function or policy. However, the convergence rate of (exact or approximate) value and policy iteration is governed by the discount factor, and so convergence is slow for discount factors close to 1 \citep{bertsekas1995,zidek2016stochastic}. Due to the particular structure of the problem we examine in our setting, we propose the use of the Binary Dynamic Search algorithm (BiDS), originally devised by \cite{mintz2017}. In contrast to value and policy iteration that operate in function spaces, the BiDS algorithm finds a fixed point in a vector space using binary search. Since the state space is small for water storage and retail management problems, our numerical results solve the exact dynamic program; however, in principle, the BiDS algorithm could also be used for approximate dynamic programming by using inexact representations of the value function.

\subsection{Outline and Contributions}
Our main contribution in this paper is to study the structure of inventory management problems with periodic review and reset control. As part of this analysis, we theoretically characterize and prove the structure of the optimal policy for these inventory management problems. Our analysis shows that the resulting policy is a threshold policy that can be seen as an extension of the classic $(s,S)$ policy. In addition, our analysis provides new theoretical methods for analyzing inventory policies since the non-convexity induced by our problem structure is distinct from that present in classical inventory models.
%

In addition to studying the structural properties of the problem, our second contribution in this paper is to employ a novel algorithm for solving stochastic optimal control problems with the specific structure of controlled resets to a single state and with constraints on the maximum time in between system resets. Unlike value and policy iteration, which require finding a fixed point in an infinite-dimensional functional space \citep{bertsekas1995}, we implement a Binary Dynamic Search (BiDS) algorithm \citep{mintz2017} that converts the control problem into finding a fixed point in a vector space using binary search. Although value and policy iteration converge exponentially, the convergence rate can be practically slow when the discount factor is close to 1 \citep{bertsekas1995,zidek2016stochastic}. The BiDS algorithm uses binary search and can thus compute the value functions of problems with controlled resets using substantially less computations. We apply BiDS to numerically solve this particular problem and experimentally verify our structural results. In doing so, we generalize the previous $(s, S)$ inventory policy to a new threshold structure that incorporates reset control. To explore the broad applicability of our new structure, we demonstrate that we have designed a policy that can be solved, interpreted, and implemented in both profit and non-profit operations.

The rest of the paper proceeds as follows. We first describe the stochastic optimization problem for the general reset control model in Section 2, with applications to water storage problem in Section \ref{sec:wat_stor_mode} and retail management problem in Section \ref{sec:ret_mngt_mode}. We then discuss the structural properties of the optimal policy and the sufficient conditions that ensure the threshold structure in Section 3, along with the numerical algorithm and its results from both water and retail problems in Section 4.

\section{Reset Control Formulation and Applications}
\label{sec:rest_control_formul}

In this section, we introduce a broad class of stochastic optimal control problems with controlled resets to a single state and with constraints on the maximum time span in between resets of the system. 
%
%
Let the subscript $n\in\mathbb{Z}_+$ denote the index of the decision epoch, and consider the discrete-time dynamical system 
\begin{equation}
\label{eqn:dyngen}
\begin{aligned}
x_{n+1} &= h(\xi_n,\tau_n,u_n,w_n)\\
t_{n+1} &= \tau_n +1\\
\xi_n &= x_n \cdot (1-r_n) + \zeta \cdot r_n\\
\tau_n &= t_n \cdot (1-r_n)
\end{aligned}
\end{equation}
where $x_n\times t_n \in \mathbb{R}^{n_x}\times\mathbb{Z}_+$ are states, $\xi_n\times \tau_n \in \mathbb{R}^{n_x}\times\mathbb{Z}_+$ are pseudo-states, $u_n\times r_n \in \mathbb{R}^{n_u}\times\mathbb{B}$ are control actions, disturbance terms $w_n \in \mathbb{R}^{n_w}$ are i.i.d. random variables, and $h : \mathbb{R}^{n_x}\times\mathbb{Z}_+\times\mathbb{R}^{n_u}\times\mathbb{R}^{n_w}\rightarrow\mathbb{R}$ is a deterministic function. The interpretation is that the control $r_n = 1$ resets the system to a \emph{known} initial state $\zeta\in\mathbb{R}^{n_x}$, the state $t_n$ keeps track of how many time steps have passed since the last system reset, and the function $h$ describes the dynamics when there is no reset. The inclusion of pseudo-states $\xi_n$ and $\tau_n$ captures the reset dynamics since these terms are set to the known initial states when $r_n = 1$.

Given the discount factor $\gamma\in[0,1)$, our goal is to solve the following stochastic control problem:
\begin{equation}
\label{eqn:genscp}
\begin{aligned}
\min\ & \mathbb{E}\Big[\textstyle \sum_{n=0}^\infty \gamma^n \Big(g(\xi_n,\tau_n,u_n,w_n) + s(x_n,t_n,w_n)\cdot r_n\Big)\Big], \\
\text{s.t. }& (\ref{eqn:dyngen}),\, t_n \leq k,\, u_n\in\mathcal{U}(\xi_n,\tau_n,w_n), \quad \text{for } n \geq0.
\end{aligned}
\end{equation}
where $g : \mathbb{R}^{n_x}\times\mathbb{Z}_+\times\mathbb{R}^{n_u}\times\mathbb{R}^{n_w}\rightarrow\mathbb{R}_+$ is a non-negative and continuous stage cost, and $s : \mathbb{R}^{n_x}\times\mathbb{Z}_+\times\mathbb{R}^{n_w}\rightarrow\mathbb{R}_+$ is a non-negative and continuous reset cost with $s(\zeta,t,w)\equiv0$.
The constraint $u_n \in \mathcal{U}(\xi_n,\tau_n,w_n)$ restricts the possible control actions to lie in a set $\mathcal{U}(\xi_n,\tau_n,w_n)$, and the constraint $t_n \leq k$ requires the system be reset at least once every $k$ time steps. For notational convenience, we will simply refer to the set $\mathcal{U}_n := \mathcal{U}(\xi_n,\tau_n,w_n)$. 

While this stochastic control formulation is quite general, we now present two instantiations of this model in real-world settings. First, we discuss a non-profit example of water storage control, and then we describe a for-profit example of retail management with changing product lines.

\subsection{Example: Water Storage Problem}
\label{sec:wat_stor_mode}

In this section, we consider the stochastic inventory control problem of water management in the developing world, where a single decision-maker must maintain the level of potable water in a residential water tank. Unlike water management in the developed world, we assume the residence does not have access to a continuous source of water, and thus it must be purchased in bulk at the beginning of the day from a communal source. Every few days, the tank needs to be fully emptied so that it can be cleaned to eliminate pathogen growth. Therefore, in this setting, the decision-maker has two actions they can take each day: first, decide if they would like to purge the tank, and second, decide how much water to purchase at the beginning of the day. These decisions must be made in such a way that optimally trades off the financial cost of purchasing water or expediting purchases when the amount is short, and the implicit health costs of letting water sit longer in the tank (that can be thought of as a holding cost).

This problem can be modeled as a stochastic control problem that is a special case of the formulation presented in Section \ref{sec:rest_control_formul}. For day $n$, let $x_n \in \mathbb{R}_+$ be the state variable representing the amount of water stored in the tank, and let $t_n \in \mathbb{Z}_+$ be the state variable representing the number of days since the tank was last emptied. Let the decision-maker's actions of the amount of water to purchase at the start of the day and whether or not to empty the tank be represented by $u_n \in \mathbb{R}_+$ and $r_n \in \mathbb{B}$, respectively. We model the demand for each day as the i.i.d. random disturbance process $w_n\in \mathbb{R}_+$.

Let the dynamics of the system be described by the following:
\begin{align}
& x_{n+1} =\big(x_n \cdot (1-r_n) + u_n - w_n\big)^+, \label{eq:x_dynam}\\
& t_{n+1} = t_n \cdot (1-r_n) + 1, \label{eq:t_dynam}
\end{align}
where $(x)^+ = \max\{x,0\}$ and $(x)^-=\min\{x,0\}$. Here, \eqref{eq:x_dynam} states that the level of water at day $n+1$ will equal the current level of water at day $n$ plus any amount of water purchased that day, minus the amount of water demanded, but cannot go below zero. If the tank is flushed on day $n$, then \eqref{eq:x_dynam} states that the water level of the tank at day $n+1$ will not depend on the amount of water at day $n$. Likewise, \eqref{eq:t_dynam} states that the number of days since the tank has been flushed increments by 1 each day until the day it is flushed, in which case the count resets to 1. Let $q: \mathbb{R}_+\rightarrow\mathbb{R}_+$ be a non-negative and non-decreasing function that represents the increased health costs associated with letting water sit in the tank for additional days. Let $c_u,c_r,p\in\mathbb{R}_{+}$ represent the variable cost of purchasing water, the per-unit cost of purging the tank, and the shortage penalty for not having sufficient water to meet demand. Using these cost parameters, state variables, dynamics equations, and discount factor $\gamma \in [0,1)$, the decision-maker's problem can be formulated as follows:
\begin{align}
\min\ & \mathbb{E}\Big[\textstyle\sum_{n=0}^\infty \gamma^n \Big(c_u u_n - p\cdot(\xi_n+u_n-w_n)^- + q(\tau_n)\cdot(\xi_n + u_n - (\xi_n+u_n-w_n)^+) + c_r x_n r_n\Big)\Big], \\
\text{s.t. }& x_{n+1} = (\xi_n + u_n - w_n)^+, \quad \text{for } n \geq 0,\\
&t_{n+1}=\tau_n + 1, \quad \text{for } n \geq 0,\\
&\xi_n = x_n \cdot (1-r_n), \quad \text{for } n \geq 0,\\
&\tau_n = t_n \cdot (1-r_n), \quad \text{for } n \geq 0,\\
&t_n \leq k, \quad \text{for } n \geq 0, \label{eq:t_n_const}\\
& r_n \in \mathbb{B}, \quad \text{for } n \geq 0,\\
&u_n \in \big[0, c_{\max} - \xi_n\big], \quad \text{for } n \geq 0. \label{eq:tank_capacity}
\end{align}
%
For this formulation, the state space is augmented to include pseudo-states $\xi_n, \tau_n$ that represent intermediate values of water in the tank and time since the last reset, given the reset action on day $n$, respectively. Using these states, the objective terms $p\cdot(\xi_n+u_n-w_n)^-$ and $q(\tau_n)\cdot(\xi_n + u_n - (\xi_n+u_n-w_n)^+)$ represent the total shortage cost and health costs at day $n$, respectively, the interpretation being that the decision-maker pays shortage costs only when demand exceeds water supply and that water quality deteriorates only if a surplus remains in the tank. The quantity $(\xi_n + u_n - (\xi_n+u_n-w_n)^+)$ is the amount of water consumed on the $n$-th day because $\xi_n + u_n$ is the amount of water available at the beginning, and $(\xi_n+u_n-w_n)^+$ is the amount of water that is still unused at the end of the day. Since the function $q$ is assumed to be monotonically increasing, the $q(\tau_n)$ term indicates that the quality of water deteriorates as time passes between emptying the tank and also ensures that consuming water that has been stored for longer durations of time is more heavily penalized. Likewise, the terms $c_u u_n$ and $c_r x_n r_n$ represent the total costs of purchasing water and flushing the tank, respectively, where the flushing cost is only incurred if the reset action is taken. The constraint \eqref{eq:t_n_const} ensures that the tank is purged at least once every $k$ days, and the constraint \eqref{eq:tank_capacity} ensures that water is not purchased in excess of the tank capacity $c_{max}$.

A defining feature of this problem is that instead of considering the holding/shortage cost trade-off, we formulate the problem based on how long water has been stored and how much water has actually been consumed. 
In circumstances where people do not have a continuous supply of water and thus have no other choice but to store water in a local storage container that lacks disinfection capabilities, we are taking an optimization approach to decide when to drain the water tank and how much water to fill it when available, with the objective of reducing the risk of contamination and shortage. In Section \ref{sec:num_res}, we will solve this problem numerically to design an efficient, interpretable policy for managing water storage systems so that it is easily implementable via a reference table that can be widely distributed to the public through paper pamphlets or the internet.

\subsection{Example: Retail Management Problem}
\label{sec:ret_mngt_mode}

Next, we consider the setting of retail inventory management in the case of changing product lines. In this setting, the decision-maker is tasked with managing the inventory of limited-time availability product lines and must decide
the quantity of stock to store and the timing for changing to a new product line such that they are able to minimize inventory costs. This is accomplished by balancing the trade-offs between the ordering cost, the holding cost for storing the inventory, the shortage cost, the stock wastage from emptying the inventory, and the cost of switching the product line.

Like the water storage problem, this retail inventory problem can also be modeled as a special case of the reset control problem presented in Section \ref{sec:rest_control_formul}. For each week $n$, let the state variables $x_n \in \mathbb{R}_+,t_n \in \mathbb{Z}_+$ represent the amount of inventory in stock of the current product line at the beginning of the week, and the number of weeks since the current product line has been offered. Let the decision-maker's weekly control actions of how much new inventory to add and whether or not to replace the current product line with a new one be given by $u_n \in \mathbb{R}_+,r_n \in \mathbb{B}$, respectively. Much like the water storage case, we model demand as an i.i.d. disturbance process $w_n \in \mathbb{R}_+$. The state dynamics for this model are identical to those in \eqref{eq:x_dynam} and \eqref{eq:t_dynam}, except that the interpretation of \eqref{eq:x_dynam} now reflects that the inventory level increases with each additional product purchased and decreases by the realization of weekly demand (or the taking of a reset action) but cannot go below zero, while \eqref{eq:t_dynam} reflects the age of the current product line.
%
Let $c_u,k_u > 0$ be the variable and fixed ordering costs associated with the current product line, respectively, and $c_r,k_r > 0$ be the associated variable and fixed ordering costs of purchasing units from the new product line. Furthermore, let $p,q>0$ be the per-unit shortage and holding costs, respectively. Then using these costs, states, dynamics, and discount factor $\gamma \in [0,1)$, the decision-maker's problem can be formulated as follows:
%
\begin{align}
    \min \ & \mathbb{E}\Big[ \sum_{n=0}^\infty \gamma^n \Big( (c_u u_n + k_u) \mathbf{1}_{\mathbb{R}^*_{+}}(u_n) - p(\xi_n+u_n-w_n)^- + q(\xi_n+u_n-w_n)^+ + (c_r x_n + k_r) r_n \Big)  \Big] \\
    \text{s.t. }& x_{n+1} = (\xi_n + u_n - w_n)^+, \quad \text{for } n \geq 0 \\
    &t_{n+1} = \tau_n + 1, \quad \text{for } n \geq 0\\
    &\xi_n = x_n \cdot (1-r_n), \quad \text{for } n \geq 0\\
    &\tau_n = t_n \cdot (1-r_n), \quad \text{for } n \geq 0\\
    &t_n \leq k, \quad \text{for } n \geq 0 \label{eq:change_const}\\
    &r_n \in \mathbb{B}, \quad \text{for } n \geq 0\\
    &u_n \in \big[0, c_{\max} - \xi_n\big], \quad \text{for } n \geq 0 \label{eq:u_caps_const}
\end{align}
where the function $\mathbf{1}_{\mathbb{R}^*_{+}}: \mathbb{R} \mapsto \mathbb{B}$ is defined as the indicator
\begin{equation}
\mathbf{1}_{\mathbb{R}^*_{+}}(u) = \begin{cases}
1, &\text{if } u > 0, \\
0, &\text{otherwise}.
\end{cases}
\end{equation}
The total cost incurred in the $n$-th week comprises up to four components: (i) the purchasing cost $c_u u_n + k_u$, where $c_u$ is the cost per unit ordered and $k_u$ is the fixed cost associated with a positive inventory order, (ii) the shortage cost $p\cdot(\xi_n+u_n-w_n)^-$ that represents the loss incurred when demand is unmet, (iii) the holding cost $q\cdot(\xi_n+u_n-w_n)^+$ for having too much inventory relative to the actual demand, and (iv) the reset cost $c_r x_n + k_r$, where $c_r$ is the waste penalty per unit discarded and $k_r$ is the fixed cost associated with resetting the product line. The constraint \eqref{eq:change_const} ensures that the product line is fully changed at least once every $k$ weeks, and the constraint \eqref{eq:u_caps_const} ensures that the inventory capacity $c_{\max}$ is not exceeded.

What differentiates this inventory management problem from the classic $(s, S)$ setting is the inclusion of reset controls. When the reset control action is taken at time $n$ (i.e., $r_n = 1$), all products are removed, and the system is reverted to the state $(x_n,t_n) = (0,0)$. This means that by solving this problem, the firm can change product lines earlier than planned if it is advantageous to do so.
%
%
In other words, given cost parameters, this problem solves the best business strategy tailored for limited-time products, which encompasses not only the optimal inventory (i.e., how much stock to order each period) but also the optimal timing (i.e., when to switch the product line). 

\section{Structural Results}
\label{sec:struct_res}

In this section, we analyze the structure of the reset control problem \eqref{eqn:genscp}. We begin by describing the dynamic programming equations first.
Let $J: \mathbb{R}^n\times \mathbb{Z}_+ \mapsto \mathbb{R}_+$  be the optimal cost-to-go function, that is, $J(x,t)$ is defined as the minimum value of \eqref{eqn:genscp} for the initial conditions $x_0 = x$ and $t_0 = t$. Let $J_0 = J(\zeta,0)$, the cost-to-go from the reset state. Then the dynamic programming equations can be characterized using the following result from \cite{mintz2017}.

\begin{proposition}
\label{prop:dpgenscp}
\citep{mintz2017} The dynamic programming equations for (\ref{eqn:genscp}) are given by
\begin{equation}
\label{eqn:ansbell}
\begin{aligned}
J(\zeta,0) &= \min_{u\in\mathcal{U}_0} \mathbb{E}\Big[ g(\zeta,0,u,w) + \gamma J(h(\zeta,0,u,w),1)\rlap{$\Big]$}\\
J(x,t) &= \min\Big\{J_0 + \mathbb{E}\big(s(x,t,w)\big), \min_{u\in\mathcal{U}_t} \mathbb{E}\Big[ g(x,t,u,w) +\gamma J(h(x,t,u,w),t+1)\Big]\Big\}\\
J(x,k) &= J_0 + \mathbb{E}\big(s(x,k,w)\big)
\end{aligned}
\end{equation}
where the middle $J(x,t)$ holds for all $x$ and $t=0,\ldots,k-1$.
\end{proposition}

These are the new dynamic programming equations that result from introducing reset control. As shown here, excluding the last period ($t=k$) at which the inventory must be reset, we evaluate every period whether resetting is more favorable. The proof for the above equations is found in \cite{mintz2017}.

Our main results in this section will prove that under a set of reasonable assumptions, the optimal fixed point policy $\pi^*:\mathbb{R}\times\mathbb{Z}_+ \mapsto \mathbb{B}\times \mathbb{R}$ to \eqref{eqn:ansbell} is a threshold policy. Specifically, we show that it can be characterized by four time-dependant parameters $s_t,S_t,\sigma_t,\Sigma_t$ and a constant $\varphi$, and has the following form: 
\begin{equation}
    \pi^*(x,t) = \begin{cases} r = 1,u=\varphi, \quad \text{ for } x\in [0,\sigma_t), \\  r=0, u = S_t - x, \quad \text{ for } x\in [\sigma_t,s_t), \\ r=0, u = 0, \quad \text{ for } x\in [s_t,\Sigma_t), \\ r=1, u = \varphi, \quad \text{ for } x\in [\Sigma_t, \infty). \end{cases}
\end{equation}
The intuition behind this policy is that in the region of $[0,\sigma_t)$, there is so little inventory that the cost of resetting is negligible, and so resetting is the optimal decision. In the region of $[\sigma_t,s_t)$, there is enough stock such that the reset cost is too high, and it is likely that the system will experience a shortage, so additional stock must be ordered. In the region $[s_t,\Sigma_t)$, there is enough stock present such that shortages are less likely than excess inventory, so no additional stock is ordered. Finally, the region $[\Sigma_t, \infty)$ indicates that the stock level is so high that it will almost surely all spoil, so it is more beneficial to reset the system.

To prove these results, we will first describe a set of sufficient assumptions that guarantee the structure of the optimal policy. Then we will present the proof for the structure of the policy, and show that it follows from the particular structure of $J$ by induction on the states.

\subsection{Technical Assumptions}

For our structural analysis of the optimal policy, we make the following technical assumptions.
\begin{assumption}
The states are $x_n\in\mathbb{R}_+$, and the state dynamics are $h(\xi_n,\tau_n,u_n,w_n) = (\xi_n+u_n-w_n)^+$.
\end{assumption}

\textbf{A1} assumes that the inventory is non-negative such that unfilled demand at each stage is not backlogged but rather lost. This translates to the system equation being $x_{n+1} = \max\{0, \xi_n + u_n - w_n\}$, instead of $x_{n+1} = \xi_n + u_n - w_n$. The state dynamics in this assumption represent the remaining inventory from the last period, which is the positive part of the stock quantity stored in the inventory plus the amount of stock added, minus the amount of stock demanded.

\begin{assumption}
The stage cost satisfies
\begin{equation}
\label{eqn:expsc}
\mathbb{E}(g(\xi_n,\tau_n,u_n,w_n)) = c u_n + K \cdot \mathbf{1}_{\mathbb{R}^*_{+}}(u_n) + H(\xi_n+u_n,\tau_n)
\end{equation}
for a function $H(\cdot,t)$ that satisfies $|\partial_x H(x+u,t)| \leq \kappa_t$ and $\partial_{xx} H(x+u,t) \geq m_t$. Here, $\partial_xH(\cdot,t)$ denotes the first derivative with respect to the first argument, and $\partial_{xx}H(\cdot,t)$ denotes the second derivative with respect to the first argument. The constants $c, K, \kappa_t \geq 0$ are non-negative, and the constants $m_t > 0$ are positive for $t = 0,\ldots,k-1$.
\end{assumption}

In \textbf{A2}, the stage cost can be decomposed into three parts, where the first term is the ordering cost, the second term is the fixed cost associated with a positive inventory order, and the third term is the function that measures the stage cost (e.g., holding/shortage cost, implying a penalty for both excess inventory and unmet demand at the end of the period). We require the function $H(\cdot,t)$ be strongly convex, and this assumption is key to proving the convexity of the value function with no reset and thereby demonstrating the optimality of the $(s, S)$ policy.

\begin{assumption}
The known initial state is $\zeta = 0$. The reset cost $R(x,t) := \mathbb{E}(s(x,t,w))$ is concave and non-decreasing in $x$. Moreover, its derivative $\partial_x R(\cdot,t)$ is bounded as $|\partial_x R(x,t)| \leq \eta_t$ and is Lipschitz continuous as
\begin{equation}
|\partial_x R(x_n,t_n) - \partial_x R(x_n',t_n)| \leq L_t\cdot|x_n-x_n'|
\end{equation}
for non-negative constants $\eta_t,L_t \geq 0$ for $t = 0,\ldots,k-1$.
\end{assumption}

\textbf{A3} assumes that we empty the entire inventory when we reset the system, and that the reset cost is concave and non-decreasing, which reflects the law of diminishing marginal utility. We also require the function $R(\cdot,t)$ be $\eta_t$-Lipshitz and $L_t$-smooth, which is key to determining the threshold that ensures the optimal policy structure.

\begin{assumption}
The input constraints are $\mathcal{U}_n = \{u : 0 \leq u \leq c_{\max} - \xi_n\}$, where the constant $c_{\max} \geq 0$ is non-negative.
\end{assumption}

\textbf{A4} is an input constraint that enforces an upper bound $c_{\max}$ on the level of stock that can be accommodated.

\begin{assumption}
The $w_n$ are i.i.d., and the density of $w_n$ is given by a function $f(w)$ that satisfies: $f(w) = 0$ for all $w < 0$, is Lipschitz continuous for $w \geq 0$ with a non-negative constant $L$, and bounded $|f(w)| \leq P$ for all $w \geq 0$ with a non-negative constant $P$.
\end{assumption}

\textbf{A5} assumes that demand is non-negative, and that its density is supported on $[0,\infty)$, Lipschitz continuous, and bounded. This assumption is fairly mild and admits a large class of probability distributions, such as exponential, gamma, and truncated normal distributions.

We will use the above five assumptions for our theoretical analysis. Furthermore, we define the following function for convenience:
\begin{equation*}
G(z, t) = c z + H(z,t) + \gamma \mathbb{E}(J((z-w)^+, t+1)).
\end{equation*}
Note that under the above definitions, we have the relation
\begin{equation*}
\min_{u \in \mathcal{U}_t} \mathbb{E}\Big[g(x,t,u,w) + \gamma J(h(x,t,u,w), t+1)\Big] = \min_{z \in [x,c_{\max}]} G(z,t) + K \cdot \mathbf{1}_{\mathbb{R}^*_{+}}(z-x) - c x.
\end{equation*}

\subsection{Proof Technique}

To prove our main result, we will use proof by induction to show that the particular structure of the policy and value function is preserved. Our approach will be to show a number of results that hold under a temporary assumption (i.e., the induction hypothesis), which will then be proved to hold under the assumptions \textbf{A1}--\textbf{A5} in a final theorem that concludes our proof. Our temporary assumption is:

\begin{assumptiont}[Assumption T.]
For any fixed $t \in \{0,\ldots,k-1\}$, we can represent $J(z, t+1)$ as
\begin{equation}
J(z, t+1) = J_i(z,t+1) \text{ if } z \in Z_{i,t} \text{ for } i = 1, \ldots, n_t,
\end{equation}
where $Z_{i,t} = [z_{i-1,t},z_{i,t}]$ forms a partition of the domain $[0,c_{\max}]$ with
\begin{equation}
0 = z_{0,t} < z_{1,t} < \cdots < z_{n,t} = c_{\max}.
\end{equation}
Moreover, the derivative of each piece $\partial_z J_i(z, t+1)$ is absolutely bounded by a finite, non-negative constant $M_{t+1}$ and is Lipschitz continuous. 
\end{assumptiont}

\textbf{AT} assumes that the value function is continuous and piecewise differentiable, and can also be partitioned into multiple pieces, the derivatives of which are Lipschitz and bounded. This assumption is required to characterize the policy structure that is complicated by the non-convexity induced by the reset control because these pieces represent the different regions of the optimal policy, which we investigate in the proceeding analysis.

\subsection{Proof of the Optimal Policy Structure}

To begin the proof, we first show a result on the smoothness of the cost-to-go function.

\begin{proposition}
\label{prop:lipsej}
If \textbf{A1}--\textbf{A5} and \textbf{AT} hold, then $\mathbb{E}(J((z-w)^+, t+1))$ has a derivative that is Lipschitz continuous with non-negative constant $(Lc_{\max}+P)M_{t+1}$, and the derivative is absolutely bounded by $M_{t+1}$.
\end{proposition}
The main proof of this proposition can be found in Appendix \ref{app:proof}. Here, we present a brief sketch. First, we show that $\mathbb{E}(J((z-w)^+,t+1))$ is differentiable using \textbf{AT} and the Leibniz integral rule for derivatives. We then show that this derivative is Lipschitz continuous by introducing an auxiliary function that has the same expectation as the cost-to-go function. Integrating over the domain and using Lemma \ref{lem:conv} gives this result. Finally, we show that the derivative is absolutely bounded by $M_{t+1}$ using standard integral and absolute value inequalities.

If the cost-to-go is piece-wise differentiable and the derivative of each piece is Lipschitz continuous, then the expected cost-to-go has a Lipschitz derivative. This is because taking the expectation smooths the function via convolution. Next, we show a result on the structure of $G(z,t)$.

\begin{proposition}
\label{prop:gztisconvex}
If \textbf{A1}--\textbf{A5} and \textbf{AT} hold, then the function $G(z,t)$ is continuous and convex on $z\in[0,c_{\max}]$ for all fixed $\gamma$ such that $0 \leq \gamma \leq m_t/((Lc_{\max}+P)M_{t+1})$.
\end{proposition}
The full proof of this proposition can be found in Appendix \ref{app:proof}, but here we present a sketch of the proof. To show this result, we note that $G$ is composed of a linear function, a twice differentiable function, and by Proposition \ref{prop:lipsej}, a function with a Lipschitz derivative. Since this implies that $G$ is absolutely continuous in its first argument, we demonstrate that it has a non-decreasing derivative, thus meaning it is convex on the desired interval.

This proposition intuitively states that the sum of the expected single stage cost (strongly convex) and the expected cost-to-go (non-convex but smooth) preserves convexity under suitable conditions, which can be computed using the parameters of the two functions. For example, $f_1(x) = x^2 /2 + \gamma \sin(x)$ is convex for $\gamma \leq 1$, whereas $f_2(x) = x^2/2 - \gamma |x|$ will never be convex unless $\gamma = 0$. This exemplifies that the value of the discount factor is critical, but the smoothness of the expected cost-to-go, which was proved in Proposition \ref{prop:lipsej}, is even more crucial in retaining the convexity of the sum.

Our next result generalizes the known result (i.e., Lemma 4.2.1 in \cite{bertsekas1995}) in inventory management models for the case of a closed feasible set.

\begin{proposition}
\label{prop:gztkcon}
If $G(\cdot,t)$ is continuous and convex on $[0,c_{\max}]$, there exists $S$ such that $G(S,t) \leq G(z,t)$ for all $z \in [0,c_{\max}]$. Furthermore, let $s = \inf \{z \in [0,S] \ | \ G(S,t) + K = G(z,t)\}$. If such $s$ exists, then 
\begin{enumerate}
\item $G(S,t) + K = G(s,t)$
\item $G(\cdot,t)$ is a non-increasing function on $[0,s]$
\item $G(S,t) + K \leq G(z,t) $ for all $z \in [0,s]$
\item $G(y,t) \leq G(z,t) + K$ for all $y,z$ with $s \leq y \leq z \leq c_{\max}$
\end{enumerate}
\end{proposition}
The full proof of this proposition can be found in Appendix \ref{app:proof}, but here we present a sketch of the proof. First, we note that because $G(\cdot,t)$ is continuous on the closed interval, it must have a minimizer, which shows the existence of $S$. Next, by continuity of $G(\cdot,t)$, we note that there may exist $s$ such that $G(s,t) = G(S,t) + K$, thus proving \textbf{C1}. We then prove \textbf{C3} using the convexity of $G(\cdot,t)$. Then using these conditions and convexity, we prove \textbf{C2}. Finally, \textbf{C4} follows from combining these results.

This proposition proves the optimality of the $(s, S)$ policy when the functions $G(\cdot,t)$ are convex and have an upper and lower bound to the allowable values of the stock. We now extend this result to the functional structure present in the reset control problem.

\begin{proposition}
\label{prop:ps}
Suppose $0 \leq \gamma \leq m_t/((Lc_{\max}+P)M_{t+1})$ for all $t = 0,\ldots,k-1$. If \textbf{A1}--\textbf{A5} and \textbf{AT} hold, then an optimal policy has the following four-stage structure with $(s_t,S_t,\sigma_t,\Sigma_t)$ thresholds:
\begin{enumerate}
\item If $x_t \in [0,\sigma_t)$, then $r_t^* = 1$ and $u_t^* = \varphi$
\item If $x_t \in [\sigma_t,s_t)$, then $r_t^* = 0$ and $u_t^* = S_t - x_t$
\item If $x_t \in [s_t,\Sigma_t)$, then $r_t^* = 0$ and $u_t^* = 0$
\item If $x_t \geq \Sigma_t$, then $r_t^* = 1$ and $u_t^* = \varphi$
\end{enumerate}
where $\varphi\in[0,c_{\max}]$ is a constant.
\end{proposition}
The full proof of this result can be found in Appendix \ref{app:proof}, but here we present a sketch. First, we define $\varphi$ as the optimal value of $u$ from the known reset state $\zeta$ with $t=0$. Because of this, from Proposition \ref{prop:dpgenscp}, we note that if the optimal reset action is $r^*(x,t)=1$, that means $u^*(x,t) = \varphi$, so we focus the proof on determining the reset policy. We now consider the set $I_R$, which is the set of all states for which a reset action would not be taken. We then consider the case when this set is empty, meaning a reset action should be taken for all states. When this set is not empty, define $\sigma_t,\sigma'_t$ as the infimum and supremum of $I_R$, respectively. Then by Proposition \ref{prop:gztkcon}, we note that there exists $S'_t = \argmin_{z\in[0,c_{max}]} G(z,t)$. However, we are not guaranteed for $s_t' = \sup \{z \in [0,S_t'] \ | \ G(S_t',t) + K \leq G(z,t)\}$ to exist. In the case where $s'_t$ does not exist, we demonstrate that the threshold $s_t$ should be set to $\sigma_t$ because $u^*(x,t)=0$ will be the minimizer for all states in this case. Next, we show that if $s'_t$ does exist, then we should set the threshold $s_t$ to $s'_t$, since then the optimal action $u^*(x,t)$ will be to order up to level $S'_t$ if $x<s'_t$, and order nothing otherwise. Next, using continuity, convexity, Proposition \ref{prop:gztisconvex}, and \textbf{A3}, we show that the region of the state space where it is optimal not to take a reset action is the interval $[\sigma_t,\sigma_t')$. This means that the optimal policy will use $\sigma_t$ value previously defined as one of the thresholds, and will use $\Sigma_t =\sigma'_t$ if $\sigma'_t < c_{max}$ and $\Sigma_t = +\infty$ otherwise. Next, we show that the threshold $S_t$ should be set to $S'_t$ by its definition and the properties of $s_t$ and $\sigma_t$.

The regions in Proposition \ref{prop:ps} represent that (i) the inventory is sufficiently small such that reset cost is negligible, and thus resetting the system and ordering a one-step quantity is optimal; (ii) the stock is sufficient such that ordering up to a given quantity and consuming the majority prior to the stock becoming too obsolete or contaminated is optimal; (iii) the stock in the system is sufficient to satisfy possible future demands, and thus resetting the system or ordering additional stock is not optimal (a.k.a., do-nothing region); and (iv) there is so much stock in the inventory that it will almost certainly never be consumed, hence making it optimal to reset the system and reorder up to a baseline amount. Next, we show how this policy influences the structure of the cost-to-go function.


\begin{proposition}
\label{prop:vf4s}
Suppose $0 \leq \gamma \leq m_t/((Lc_{\max}+P)M_{t+1})$ for all $t = 0,\ldots,k-1$. If \textbf{A1}--\textbf{A5} and \textbf{AT} hold, then the function $J(x,t)$ for $t = 0,\ldots,k-1$ has the following form with $(s_t,S_t,\sigma_t,\Sigma_t)$ thresholds:
\begin{equation}
\label{eqn:jtz}
J(x,t) = \begin{cases}
J_0 + \mathbb{E}(s(x,t,w)), & x \in [0,\sigma_t) \\
c (S_t-x) + K + H(S_t,t) + \gamma \mathbb{E}(J((S_t-w)^+,t+1)), & x \in [\sigma_t,s_t) \\
H(x,t) + \gamma \mathbb{E}(J((x-w)^+,t+1)), & x \in [s_t,\Sigma_t) \\
J_0 + \mathbb{E}(s(x,t,w)), & x \geq \Sigma_t
\end{cases}
\end{equation}
Furthermore, each piece has a Lipschitz derivative, and the derivative of each piece is absolutely bounded by $M_t = \eta_t + c + \kappa_t + m_t/(Lc_{\max} + P)$ for $t = 0,\ldots,k-1$.
\end{proposition}
The full proof of this proposition can be found in Appendix \ref{app:proof}. Here, we present a proof sketch. The functional structure of the policy follows directly as a consequence of Proposition \ref{prop:ps}, so we focus on showing the properties of the derivative of the cost-to-go function. We thus examine each piece of the function individually. The result is proved for the first and fourth pieces as a consequence of \textbf{A3}, and we note that in the second region, $J(x,t)$ is linear and thus also satisfies the result. For the third piece, the result follows as a consequence of Proposition \ref{prop:lipsej} and \textbf{A2}.

\begin{figure}[h]
\centering
{\renewcommand{\arraystretch}{0.7}
\begin{tabular}{p{3.9cm} p{3.7cm} p{3.7cm} p{3.7cm}}
\multicolumn{4}{c}{\includegraphics[width=16cm]{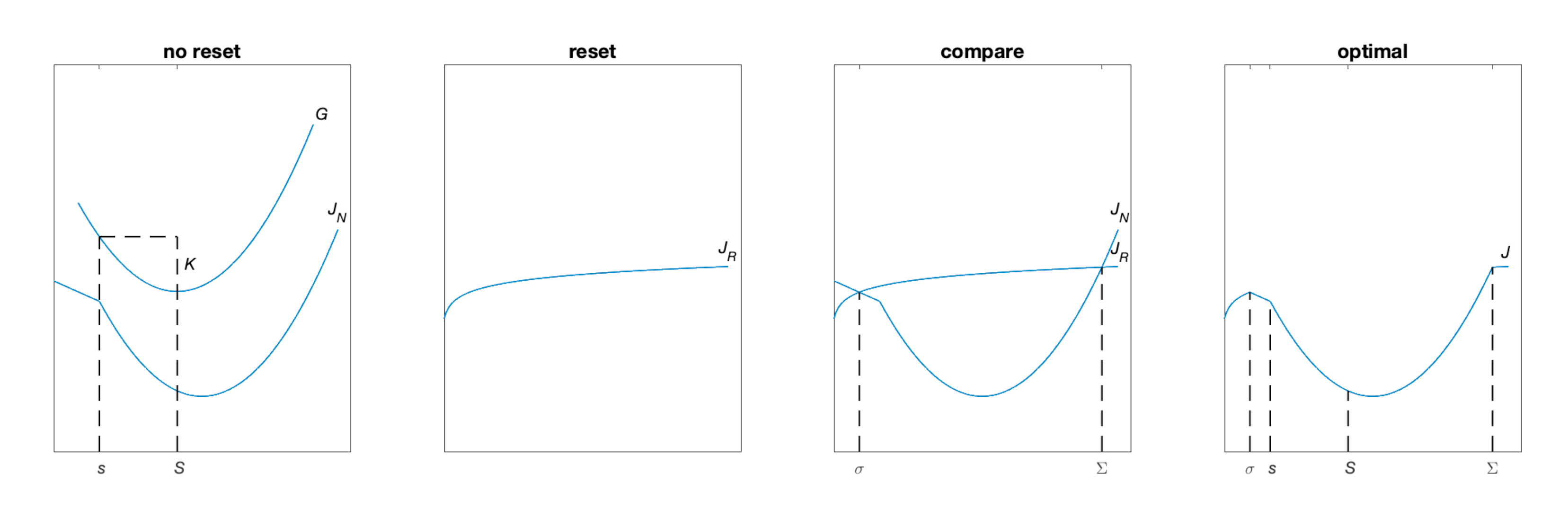}} \\
\hline
\hfil linear then convex & \hfil concave & \hfil points of intersection & \hfil continuous, piece- \\
\hfil $(s, S)$ thresholds &  & \hfil $(\sigma, \Sigma)$ thresholds & \hfil wise differentiable \\
\hline \\
\end{tabular}
}
\caption{Proof sketch for Proposition \ref{prop:vf4s}}
\label{fig:vf4s}
\end{figure}

The policy structure from Proposition \ref{prop:ps} implies that the value function has up to four thresholds. The first two thresholds, namely $s$ and $S$, are obtained from the value function when there is no reset, which is linear and then convex. The points of intersection that result from juxtaposing the value functions with and without reset give us the last two thresholds, $\sigma$ and $\Sigma$, thereby yielding the optimal cost-to-go that is continuous and piecewise differentiable, as shown in Figure \ref{fig:vf4s}. We now complete the proof by induction and show that the optimal policy and cost-to-go functions indeed follow the forms above.

\begin{theorem}
\label{thm:fin}
Suppose $\gamma$ is such that $0 \leq \gamma \leq \min \{\gamma_t\ |\ t\in\{0,\ldots,k-1\}\}$ for
\begin{equation}
\begin{aligned}
\gamma_{k-1} &= \frac{m_{k-1}}{(Lc_{\max}+P)\cdot \eta_k}\\
\gamma_t &= \frac{m_t}{(Lc_{max}+P)\cdot(\eta_{t+1} + c + \kappa_{t+1}) + m_{t+1}}\ \mathrm{for}\ t = 0,\ldots,k-2
\end{aligned}
\end{equation}
If \textbf{A1}--\textbf{A5} hold, then the policy described in Proposition \ref{prop:ps} is optimal and the value function has the structure described in Proposition \ref{prop:vf4s}.
\end{theorem}
The complete proof of the theorem can be found in Appendix \ref{app:proof}, but here we present a sketch. First, we note that when $t=k$, the optimal policy will choose to reset the system for all $x \in [0, c_{max}]$, which means $J(x,k) = J_0 + \mathbb{E}(s(x,k,w))$. The rest follows by Propositions \ref{prop:ps} and \ref{prop:vf4s} and induction on $t$.

This concludes our proof by induction, from which we obtain the threshold for the discount factor that guarantees the structural properties of the optimal policy that we proved using assumptions \textbf{A1}-\textbf{A5}. We can also safely remove the temporary assumption \textbf{AT}, as Proposition \ref{prop:vf4s} proves that the value function can be partitioned into at most four pieces, each of which is Lipschitz continuous and absolutely bounded.

\subsection{Policy Structure for Water Storage and Retail Management Problems}

In this section, we provide sufficient conditions for both water storage and retail problems to show when their policies have the four-threshold structure of the reset problem. To do this, we show what conditions are sufficient such that each model satisfies assumptions \textbf{A1}-\textbf{A5}. First, we outline the conditions under which the water storage problem satisfies all assumptions.

\begin{proposition}
\label{prop:sctc}
Suppose \textbf{A5} (the assumption about the distribution of $w_n$) holds, and that $f(w) > 0$ for all $w\in[0,c_{\max}]$. If $p-q(t) > 0$ for all $t = 0,\ldots,k-1$, then the water storage problem described in Section \ref{sec:wat_stor_mode} satisfies \textbf{A1}--\textbf{A5}.
\end{proposition}
The main proof of this proposition can be found in Appendix \ref{app:proof}, but here we present a sketch. First, we note that \textbf{A1} and \textbf{A4} hold by the definition of the problem, and that \textbf{A5} holds by the assumptions of the proposition. Then we show that through a reformulation of the water problem stage cost (and setting the fixed cost to zero), the problem can be shown in conjunction with the Leibniz integral rule to satisfy \textbf{A2}. Next, we show that \textbf{A3} follows from the structure of the reset conditions in this problem.

The following proposition presents sufficient conditions for the retail problem to have the same threshold policy structure.

\begin{proposition}
\label{prop:sctc2}
Suppose \textbf{A5} holds, and that $f(w) > 0$ for all $w\in[0,c_{\max}]$. Then the retail management problem described in Section \ref{sec:ret_mngt_mode} satisfies \textbf{A1}--\textbf{A5}.
\end{proposition}
The main proof of this proposition can be found in Appendix \ref{app:proof}. However, we present a sketch here. First, we note that \textbf{A1}, \textbf{A4}, and \textbf{A5} hold by definition of the problem and by assumption of the proposition. Next, letting the fixed cost $K=k_u$ and using the Leibniz integral rule, the problem can be shown to satisfy \textbf{A2}. Finally, \textbf{A3} follows from the reset cost structure.

\section{Numerical Results}
\label{sec:num_res}

In this section, we present numerical studies to validate the theoretical structure we derived in Section \ref{sec:struct_res}. All computations are performed in MATLAB 2018b on a laptop computer with a 2.6GHz processor and 16GB of RAM. First, we present the numerical dynamic programming algorithm that we will use to compute the optimal value functions and policies. We then explore the two case studies of water storage and retail management problems and numerically verify that their value functions and optimal policies follow the theoretical structure.
\subsection{Binary Dynamic Search Algorithm}
\begin{algorithm}[h]
\begin{algorithmic}[1] 
\caption{Binary Dynamic Search (BiDS) Algorithm \citep{mintz2017}}
\label{alg:bids_alg}
\State initialize $ \underline{v} \leftarrow 0$ and $\overline{v} \leftarrow (\frac{1}{1-\gamma})\min_u \mathbb{E}[g(\zeta,0,u,w_0) + 
\gamma\cdot s(h(\zeta,0,u,w_0),1,w_1)]$
\Repeat
\State set $v \leftarrow (\overline{v}+\underline{v})/2$
\State set $V(x,k,v) = v + \mathbb{E}[s(x,w,k)]$
\For {$t = (k-1),(k-2),\ldots,0$}
\State set $V(x,t,v) = \min\big\{v + \mathbb{E}[s(x,w,t)], \min_{u\in{\mathcal{U}_t}} \mathbb{E}[g(x,t,u,w) + \gamma V(h(x,t,u,w),t+1,v)]\big\}$
\EndFor
\State set $\Upsilon(v) = \min_{u\in{\mathcal{U}_0}} \mathbb{E}[g(\zeta,0,u,w) + \gamma V(h(\zeta,0,u,w),1,v)]$
\If{$v > \Upsilon(v)$}
\State set $\overline{v}\leftarrow v$
\Else
\State set $\underline{v}\leftarrow v$
\EndIf
\Until {$(\overline{v}-\underline{v}) \leq \epsilon$}
\State set $v^*=(\overline{v}+\underline{v})/2$
\end{algorithmic}
\end{algorithm}

We solve the dynamic programming equations and compute the optimal policy using the Binary Dynamic Search (BiDS) algorithm, initially developed by \cite{mintz2017}. The classical discounted reward settings, optimal value functions, and policies are computed with some form of value iteration (VI) and policy iteration (PI). However, if the state space is infinite, this computation is equivalent to finding a fixed point in an infinite-dimensional functional space, and while strong convergence guarantees exist in a discounted setting, the rate of convergence is highly dependent on the discount factor and can be practically slow \citep{bertsekas1995}. For these reasons, these problems are generally solved using approximate dynamic programming methods in practice. However, the BiDS algorithm can exploit the specific structure of reset control problems to calculate the optimal value function and policy with arbitrary precision. To do this, BiDS converts the problem into finding a fixed point in a vector space using binary search, where this vector can be thought of as the optimal value function evaluated at the reset state and time 0. As can be seen in Algorithm \ref{alg:bids_alg}, BiDS first initializes the search space using an upper and lower bound on the value function at the reset state. Then BiDS takes the midpoint of the interval and uses backward induction to compute the implied reset state cost. Depending on how this value compares with the candidate for the iteration, a new search interval is selected using the same procedures as a binary search. The algorithm terminates either when a true fixed point is found or when the numerical tolerance $\epsilon$ is reached. The theoretical convergence and computational guarantees of BiDS can be found in \cite{mintz2017}.

\subsection{Example: Water Storage Problem}

We use BiDS to numerically solve the water storage problem from Section \ref{sec:wat_stor_mode}. We generate random demand from the truncated normal distribution that is bounded from below by zero \citep{botev2017} and set $k = 7$ to ensure that we flush at least once a week. The cost parameters are carefully chosen such that we place more weight on shortage cost, compared to purchasing and flushing costs (i.e., $p > c_u > c_r$). Furthermore, the cost of consuming water that has been stored long-term exceeds the cost of flushing water (i.e., $q(t) < p$ for all $t$, but $q(t) \geq c_r$ for $t$ close to $k$) such that the model favors resetting the tank more often to prevent contamination. In addition, while exponential microbial growth is normally assumed \citep{maier2009}, we present a simplified model with linear growth to account for fluctuations in the microbial load that result from water usage and refilling. Note that any other health penalty could be used as long as it satisfies the monotonicity condition.

\begin{figure}[h]
\centering
\includegraphics[width=0.6\textwidth]{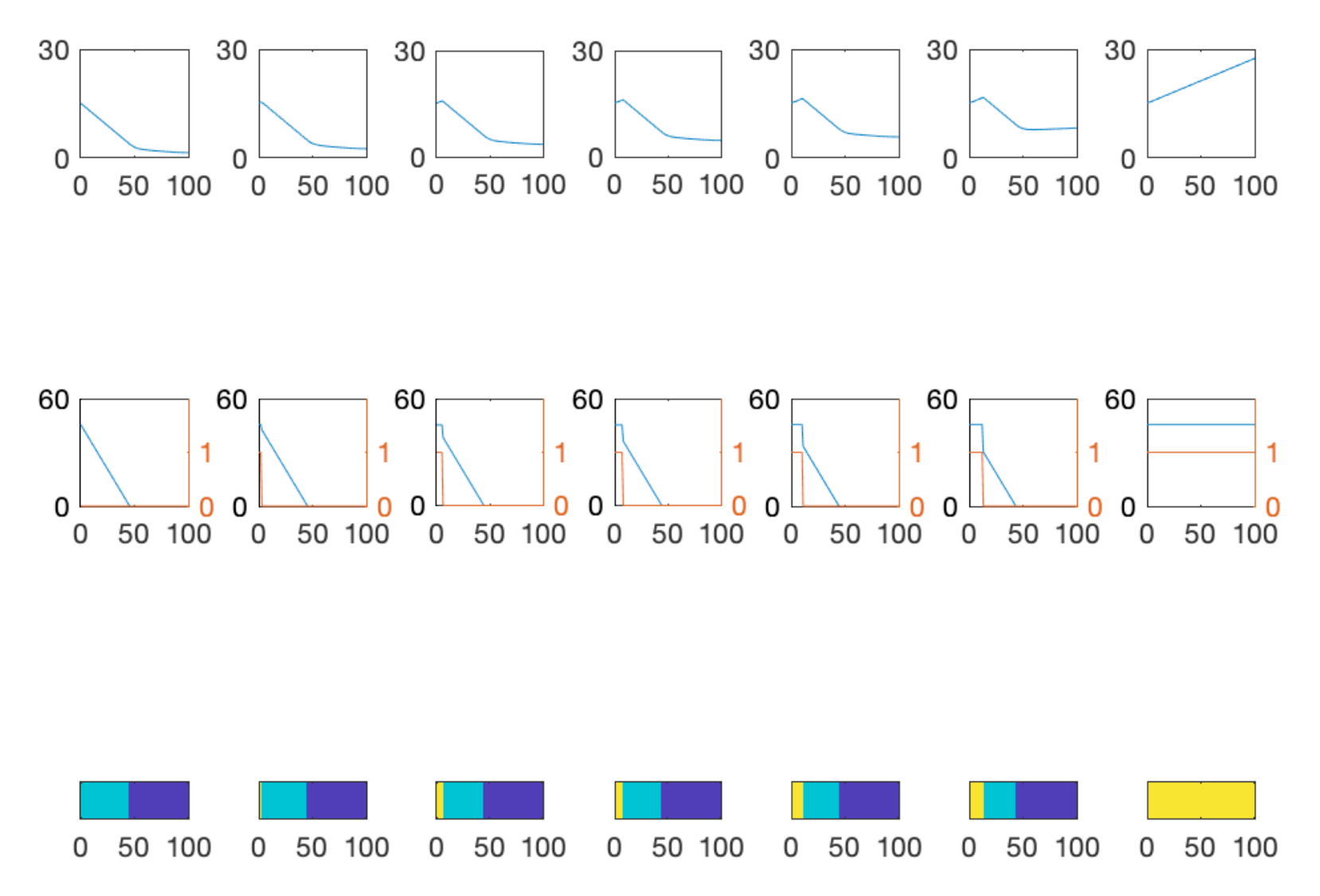}
\caption{Example of an optimal policy for the water storage problem}
\label{fig:water_policy_1}
\end{figure}

In Figure \ref{fig:water_policy_1}, the top row shows the value function, and the second row shows the optimal control actions in blue (optimal fill amount) and red (reset control action). The bottom row shows the different zones for the optimal policy---yellow corresponds to flushing the tank and reordering water, cyan to ordering water and not flushing, and dark blue to not doing anything. For all plots, the state $x$ (the amount of water in the tank) is on the $x$-axis, and the subplots from left to right correspond to the states $t = 1, \ldots, 7$ (the number of days since the tank was last emptied).

\begin{figure}[b]
\centering
\begin{subfigure}[b]{0.6\textwidth}
\includegraphics[width=1\linewidth]{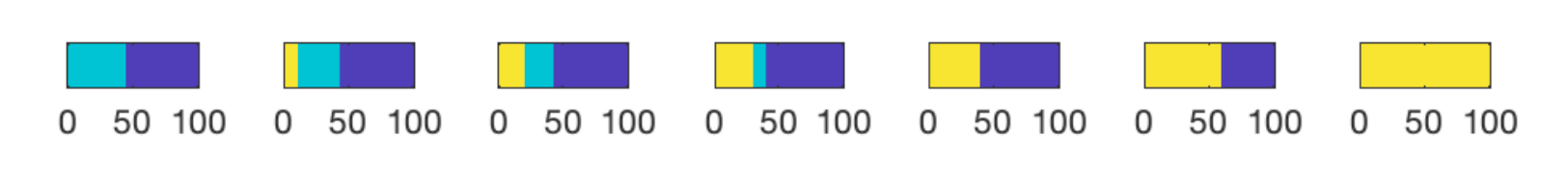}
\label{fig:ex2}
\end{subfigure}
\begin{subfigure}[b]{0.6\textwidth}
\includegraphics[width=1\linewidth]{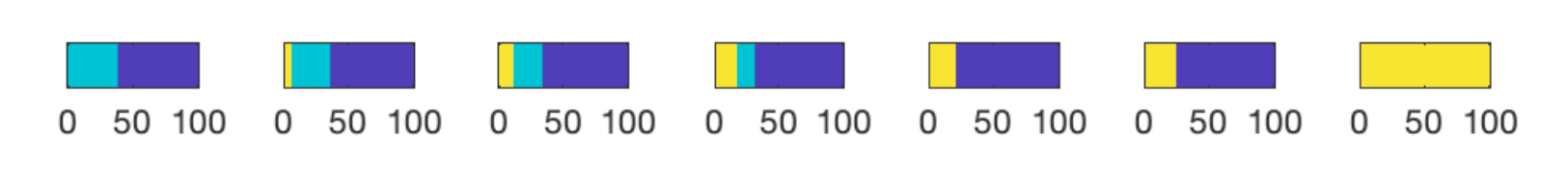}
\label{fig:ex3}
\end{subfigure}
\caption{Optimal policy when there is a high risk of contamination (top), exacerbated by limited access to water (bottom)}
\label{fig:water_policy_2}
\end{figure}

With all other parameters fixed from Figure \ref{fig:water_policy_1}, we analyze the following two scenarios in Figure \ref{fig:water_policy_2}. When there is a high risk of contamination (top), we impose a larger penalty on water consumption to reflect this condition. We can observe that the optimal policy indeed suggests flushing the tank more often, as indicated by thicker yellow and thinner cyan bands, to reduce the risk of consuming contaminated water. Further exacerbation due to limited access to water (bottom), which is often the case in the developing world \citep{ohchr2015, wateraid2016}, can be modeled by assigning a higher cost to purchasing water. Then we can see that the control policy is mostly composed of reset (yellow) and do-nothing (dark blue) regions, with their thresholds slightly shifted to the left compared to the previous policy. As illustrated by wider blue and narrower yellow regions, this regrettably implies that we cannot afford to reset that often.
Nevertheless, our model provides a guideline that we can rely on, especially when we are faced with a difficult choice, whether it suggest a breakthrough or a compromise.

\subsection{Example: Retail Management Problem}

\begin{figure}[h]
\centering
\begin{subfigure}[b]{0.6\textwidth}
\includegraphics[width=1\linewidth]{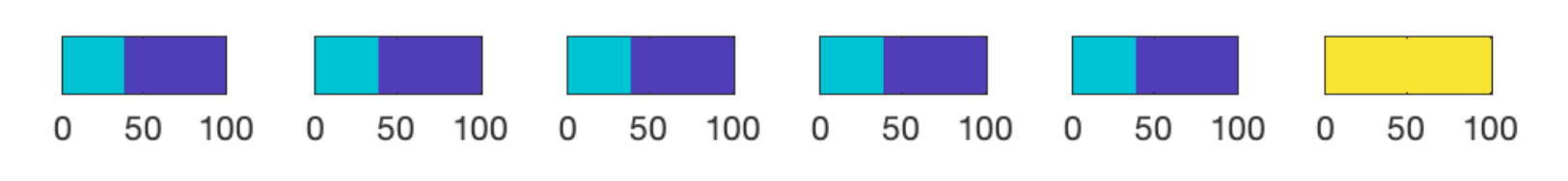}
\label{fig:ex4}
\end{subfigure}
\begin{subfigure}[b]{0.6\textwidth}
\includegraphics[width=1\linewidth]{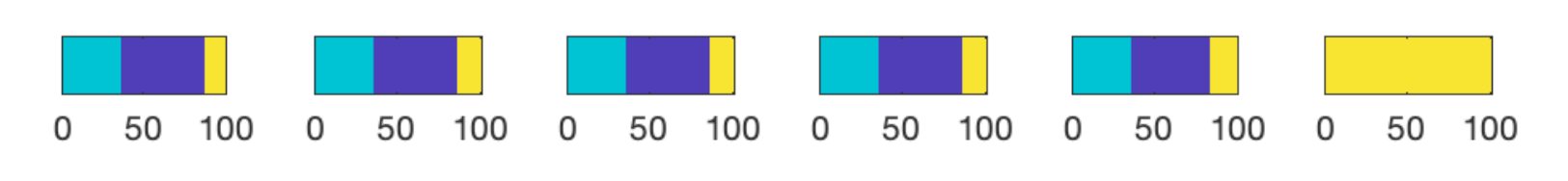}
\label{fig:ex5}
\end{subfigure}
\caption{Optimal policy for different cases of the retail management problem}
\label{fig:retail_policy_1}
\end{figure}

In this section, we solve the retail management problem described in Section \ref{sec:ret_mngt_mode} using the BiDS algorithm under two different cases, as shown in Figure \ref{fig:retail_policy_1}. If a firm lacks the flexibility required to swiftly change the product line (top), we model this by assigning a higher reset cost. We can then observe that the optimal policy essentially reduces to an $(s, S)$ policy, similar to that of classical inventory management problems, where firms do not have the infrastructure and resources to switch the product line often. In contrast, if a company is situated in a rapidly evolving industry and thus equipped with the speed and adaptability demanded by the market (bottom), we use a lower reset cost but a higher holding cost to take the ephemeral nature of our merchandise into consideration. We can now see that the inventory policy has a yellow region on the right. This encourages us to change inventory when there is too much on hand, and therefore epitomizes the right response to the transiency of a trend. 

At each $t$, there are at most three thresholds that separate the control policy into four regions. In the first region (yellow), the trivial amount of stock renders the inventory reset cost negligible. Hence, emptying the inventory and ordering an optimal one-step quantity is optimal. In the next region (cyan), the sufficiency of stock in the inventory makes it optimal to order up to some quantity and still consume most of the stock before it becomes obsolete. In the third `do-nothing region’ (blue), it is not optimal to reset the inventory or order further stock as potential future demand can already be satisfied. In the final region (yellow), the inventory stock is almost certain to never be consumed, and it is thus optimal to empty the inventory and reorder up to some optimal quantity, changing the product line if appropriate. Also, note that this covers all possible stages that can be suggested by the optimal policy, and that the actual policy may not have all of them, as in the case of Figure \ref{fig:retail_policy_1}.

\begin{figure}[h]
\centering
\includegraphics[width=0.6\textwidth]{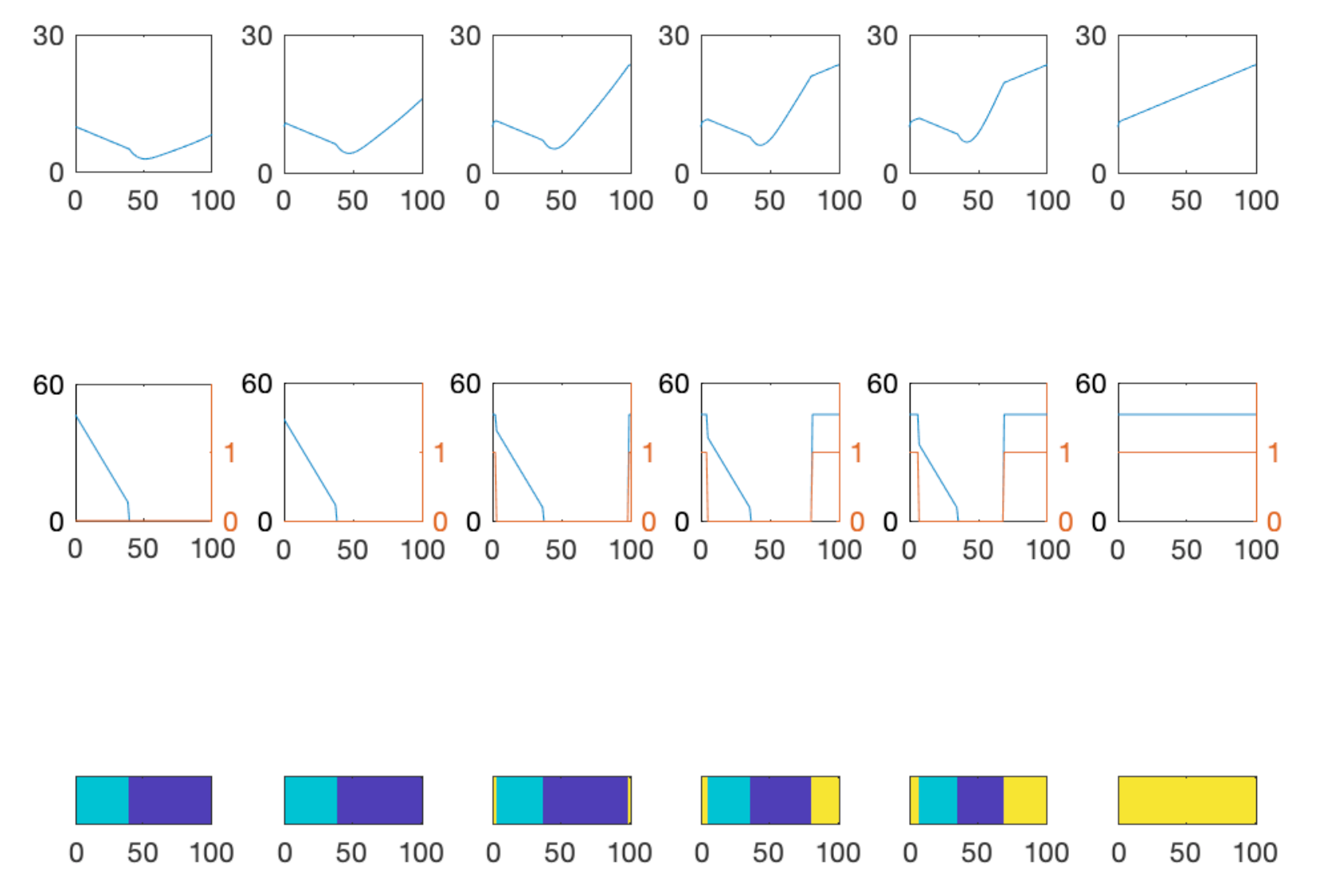}
\caption{Example of an optimal policy for deteriorating item with linear time-dependent holding cost}
\label{fig:retail_policy_2}
\end{figure}

\section{Conclusion}

Inventory management that empowers us to flexibly adapt to change is necessary to survive in a market and society that continues to evolve with growing acceleration. Previous inventory policies were contingent upon the tradeoff balance between holding and shortage costs, and informed us only of the optimal inventory position, or the quantity that must be ordered, for each period. By introducing the option to reset inventory into such classical frameworks, we are pioneering the design of more dynamic policies that overcome the limitations of pre-existing static policies. Furthermore, this enables us to continuously assess whether resets are more profitable at the preset periodic intervals or at earlier time points. In this paper, we (i) theoretically investigated the structural properties underlying the optimal policy of this problem, and (ii) implemented the appropriate algorithm to numerically solve and empirically validate the policy structure of the aforementioned problem. Finally, we demonstrated the broad utility of this problem by providing both non-profit and for-profit examples.

The original motivation for this paper stemmed from the local water inventory management for developing nations \citep{mintz2017} as they often lack water distribution networks that provide clean water throughout the day. As a result, homes and apartments utilize water storage systems that are filled during a small window of time in the day when the water distribution network is active. However, these water storage systems are not equipped with disinfection capabilities. This leads to water being stored for long durations, thereby exposing people to substantially increased bacterial and viral water contamination \citep{coelho2003, lee2005, kumpel2016}. Furthermore, obtaining even the minimal quantity of water necessary for survival can be a huge burden for people without a continuous water supply at home. This forces these individuals to compromise their health by using less water or collecting water from unsafe sources \citep{wateraid2016}. Our model seeks to promote universal and equitable access to safe water by providing guidance towards making decisions when confronted with such healthcare dilemmas.

Our formulation and proof are also sufficiently flexible to account for general reset control problems in supply chain management. For instance, we can incorporate modifications, such as introducing a time-varying holding cost to take the depreciation rate into account and penalize based on the amount of time it has been sitting on a shelf, with the purpose of managing aging inventories (see Figure \ref{fig:retail_policy_2}). Finally, our results suggest that the sufficient conditions that ensure the threshold structure of the optimal policy are marginally conservative since they are not necessary conditions, and that there may exist a relaxation of these conditions, which we leave for future studies.

\ACKNOWLEDGMENT{This material is based upon work supported by the National Science Foundation under Grant CMMI-1847666.
}

%
%
%

\begin{APPENDICES}

\numberwithin{equation}{section}

\section{Proof of Propositions and Theorem}

\label{app:proof}
\proof{Proof of Proposition \ref{prop:lipsej}: }
We first show that $\mathbb{E}(J((z-w)^+,t+1))$ is differentiable. Consider any $z \in [0,c_{\max}]$, and observe that
\begin{equation}
\begin{aligned}
J((z-w)^+,t+1) &= \begin{cases} 
J_1(0,t+1), & \text{if } z-w \leq 0 \\
J_i(z-w, t+1), & \text{if } z-w \in Z_{i,t}
\end{cases} \\
&= \begin{cases} 
J_1(0, t+1), & \text{if } w \geq z \\
J_i(z-w, t+1), & \text{if } w \in [z-z_{i,t}, z-z_{i-1,t}] 
\end{cases}
\end{aligned}
\end{equation}
Now, let $m$ be the smallest integer such that $z \in Z_{m,t}$. Then the expectation is given by
\begin{multline}
\mathbb{E}(J((z-w)^+,t+1)) = \int_0^\infty J((z-w)^+,t+1) f(w) dw = \int_0^{z-z_{m-1,t}} J_{m}(z-w,t+1) f(w) dw + \\
\sum_{i=1}^{m-1} \int_{z-z_{i,t}}^{z-z_{i-1,t}} J_{i}(z-w,t+1) f(w) dw + \int_z^\infty J_{1}(0,t+1) f(w) dw.
\end{multline}
Since by assumption both $J_{i}(z,t+1)$ and $\partial_z J_{i}(z,t+1)$ are continuous for $z \in Z_{i,t}$, using the Leibniz integral rule gives
\begin{multline}
\partial_z \mathbb{E}(J((z-w)^+,t+1)) = J_{m} (z_{m-1,t},t+1) f(z-z_{m-1,t}) + \int_0^{z-z_{m-1,t}} \partial_z J_{m}(z-w,t+1) f(w) dw + \\
\sum_{i=1}^{m-1} \Big(J_{i}(z_{i-1,t},t+1) f(z-z_{i-1,t+1}) - J_{i}(z_{i,t},t+1) f(z-z_{i,t}) + \int_{z-z_{i,t}}^{z-z_{i-1,t}} \partial_z J_{i}(z-w, t+1) f(w) dw \Big) - \\
J_{1}(0,t+1) f(z) + \int_z^\infty \partial_z J_{1}(0,t+1) f(w) dd = \int_0^{z-z_{m-1,t}} \partial_z J_{m}(z-w,t+1) f(w) dw + \\
\sum_{i=1}^{m-1} \int_{z-z_{i,t}}^{z-z_{i-1,t}} \partial_z J_{i}(z-w,t+1) f(w) dw + \int_z^\infty \partial_z J_{1}(0,t+1) f(w) dw = \\
\int_0^\infty \partial_z J((z-w)^+, t+1) f(w) dw
\end{multline}
where the last equality follows from the fact that $f(\cdot)$ is Lipschitz continuous and $J(\cdot,t+1)$ is differentiable almost everywhere on its domain since it is piecewise differentiable by assumption. Hence, we can conclude that $\mathbb{E}(J((z-w)^+,t+1))$ is differentiable.

Next, we show that this derivative is Lipschitz continuous. Define
\begin{equation}
J_A(x,t+1) = \begin{cases} 
J(0,t+1), & x \leq 0 \\
J(x,t+1), & x \in [0,c_{\max}] \\
J(c_{\max},t+1), & x \geq c_{\max}
\end{cases}
\end{equation}
where the subscript $A$ indicates ``auxiliary''. Then for $z \in [0,c_{\max}]$ and $w \in [0,\infty)$, we have $J((z-w)^+,t+1) = J_A(z-w,t+1)$. Note that $\mathbb{E}(J_A(z-w,t+1)) = \int_0^\infty J_A(z-w,t+1) f(w) dw = \int_{-\infty}^\infty J_A(z-w,t+1) f(w) dw = (J_A(\cdot,t)*f)(z)$ (i.e., the convolution of $J_A(\cdot,t+1)$ and $f(\cdot)$), where the second equality follows since $f(w) = 0$ for $w \in (-\infty,0)$ by assumption. This implies that
\begin{equation}
\label{eqn:conv}
\begin{split}
\partial_z (J_A(\cdot,t+1)*f)(z) &= \partial_z \mathbb{E}(J_A(z-w,t+1)) \\
&= \partial_z \mathbb{E}(J((z-w)^+,t+1)) \\
&= \int_0^\infty \partial_z J((z-w)^+,t+1) f(w) dw \\
&= \int_0^\infty \partial_z J_A(z-w,t+1) f(w) dw \\
&= (\partial_z J_A(\cdot,t+1) * f)(z)
\end{split}
\end{equation}
Observe that 
\begin{equation}
\partial_z J_A(z,t+1) = \begin{cases}
0, & x < 0 \\
\partial_z J(z,t+1), & z \in (0,c_{\max}) \\
0, & x > c_{\max}
\end{cases}
\end{equation}
This means $\int_{-\infty}^\infty|\partial_z J_A(z,t+1)|dz = \int_0^{c_{\max}}|\partial_zJ(z,t+1)|dz \leq c_{\max}M_{t+1}$ since we had assumed $|\partial_z J(z,t+1)| \leq M_{t+1}$. Thus, Lemma \ref{lem:conv} implies $\partial_z \mathbb{E}(J((z-w)^+,t+1))$ is Lipschitz with constant $(Lc_{\max}+P)M_{t+1}$.

Lastly, we show that the derivative is absolutely bounded by $M_{t+1}$. From (\ref{eqn:conv}), we have
\begin{equation}
\begin{aligned}
|\partial_z (J_A(\cdot,t+1)*f)(z)| & = |(\partial_x J_A(x,t+1) * f)(z)|\\
&= \Bigg|\int_0^\infty \partial_z J((z-w)^+,t+1) f(w) dw\Bigg| \\
&\leq \int_0^\infty \big|\partial_z J((z-w)^+,t+1)\big|\cdot \big|f(w)\big| dw \\
&\leq M_{t+1}\int_0^\infty \big|f(w)\big| dw \\
&\leq M_{t+1}
\end{aligned}
\end{equation}
where in the last line we have used the facts that $f(w)$ is a probability density, and hence non-negative, and integrates to one. $\Halmos$
\endproof

\proof{Proof of Proposition \ref{prop:gztisconvex}: }
We first recall the definition
\begin{equation*}
G(z, t) = c\cdot z + H(z, t) + \gamma \mathbb{E}(J((z-w)^+, t+1)).
\end{equation*}
Note that $c\cdot z$ is linear, that $H(z,t)$ is twice differentiable by assumption, and that $\mathbb{E}(J((z-w)^+, t+1))$ has a Lipschitz derivative by Proposition 2. Hence, $G(\cdot,t)$ is absolutely continuous in its first argument. Since the domain $[0,c_{\max}]$ is closed and bounded, this means that $G(\cdot,t)$ is convex if it has a non-decreasing derivative (see page 115 of \cite{fremlin2001}). From \textbf{A2}, we have that $\partial_{zz}(c\cdot z + H(z,t)) \geq m_t$, and Proposition 2 gives that $\partial_z \mathbb{E}(J((z-w)^+,t+1))$ is Lipschitz with constant $(Lc_{\max}+P)M_{t+1}$. Thus, by Lemma \ref{lemma:nodec}, we have that $\partial_z G(z,t)$ is non-decreasing for $\gamma$ such that $0 \leq \gamma \leq m_t/((Lc_{\max}+P)M_{t+1})$. $\Halmos$
\endproof

\proof{Proof of Proposition \ref{prop:gztkcon}: }
Since $G(\cdot,t)$ is continuous on the closed interval $[0,c_{\max}]$, there exists a minimizer of $G(\cdot,t)$. Let $S \in \arg \min_{z \in [0,c_{\max}]} G(z,t)$. For the remainder of the proof, we consider the case where $s$ exists. By the continuity of $G(\cdot,t)$ and the definition of $s$, we must have $G(S,t) + K = G(s,t)$. Next, for all $z \in [0,s]$, there exists $\lambda \in [0,1]$ such that $s = \lambda z + (1-\lambda) S$ and
\begin{equation}
\label{eqn:gstleqgzt}
G(s,t) \leq \lambda G(z,t) + (1-\lambda) G(S,t) \leq \lambda G(z,t) + (1-\lambda) G(z,t) = G(z,t),
\end{equation}
where the inequalities follow from the convexity of $G(\cdot,t)$ and the definition of $S$, respectively. Next, consider any $z_1$ and $z_2$ with $0 \leq z_1 \leq z_2 \leq s$, and note that there exists $\lambda \in [0,1]$ such that $z_2 = \lambda z_1 + (1-\lambda) s$. Then we have
\begin{equation}
G(z_2,t) \leq \lambda G(z_1,t) + (1-\lambda) G(s,t) \leq \lambda G(z_1,t) + (1-\lambda) G(z_1,t) = G(z_1,t),
\end{equation}
where the second inequality follows by (\ref{eqn:gstleqgzt}). This shows $G(\cdot,t)$ is non-increasing on $[0,s]$. Since $G(\cdot,t)$ is continuous, this means by definition of $s$ that $G(S,t) + K \leq G(s,t)$. Using the above result that $G(\cdot,t)$ is non-increasing on $[0,s]$, we have $G_t(S) + K \leq G_t(z)$ for all $z \in [0,s]$. By definition of $s$, for any $y \in [s,c_{\max}]$, we have $G_t(S) + K \geq G_t(y)$. However, we also know that $G_t(S) \leq G_t(z)$ for any $z \in [y,c_{\max}]$. Combining the two yields the last part of the result. $\Halmos$
\endproof

\proof{Proof of Proposition \ref{prop:ps}: }
We first define the constant
\begin{equation}
\varphi = \arg\min_{u\in\mathcal{U}_0} \mathbb{E}\Big[g(\zeta,0,u,w) + \gamma J(h(\zeta,0,u,w),1)\Big].
\end{equation}
From Proposition 1, we have $u_t^*(x) = \varphi$ whenever $r_t^*(x) = 1$. So, our proof will be structured around determining an optimal reset policy $r_t^*(x)$. We define the functions
\begin{equation}
\begin{aligned}
J_R(x,t) &= J_0 + \mathbb{E}(s(x,t,w))\\
J_N(x,t) &= \min_{u \in [0,c_{\max}-x]} G(x+u,t) + K \cdot \mathbf{1}_{\mathbb{R}^*_{+}}(u) - c\cdot x.
\end{aligned}
\end{equation}
where the subscript $R$ indicates ``reset'', and the subscript $N$ indicates ``no reset''. Observe that these are the value functions corresponding to $r_t = 1$ and $r_t = 0$, respectively. Next, define the set $I_R = \{x \in [0,c_{\max}] \ | \ J_R(x,t) \geq J_{N}(x,t)\}$, and observe that $I_R$ is bounded by construction. Now, we consider two cases:

The first case is when $I_R = \emptyset$. Then by definition of $I_R$, we have $J_R(x,t) < J_{N}(x,t)$ for all $x \in [0,c_{\max}]$. For this case, an optimal policy is to choose $\sigma_t = s_t = \Sigma_t = +\infty$ since it is optimal to choose $r_t^*(x) = 1$ for all $x\in[0,c_{\max}]$.

The second case is when $I_R \neq \emptyset$. Let $\sigma_t = \inf \{x \in [0,c_{\max}] \ | \ J_R(x,t) \geq J_{N}(x,t)\}$ and $\sigma_t' = \sup \{x \in [0,c_{\max}] \ | \ J_R(x,t) \geq J_{N}(x,t)\}$, and note that $\sigma_t$ and $\sigma_t'$ are finite since $I_R$ is bounded and non-empty. Let $S_t' = \argmin_{z\in[0,c_{\max}]}G(z,t)$ and define $s_t' = \sup \{z \in [0,S_t'] \ | \ G(S_t',t) + K \leq G(z,t)\}$. Note that by Proposition 4, $S_t'$ is guaranteed to exist whereas $s_t'$ may or may not exist. We consider two subcases based on the existence of $s_t'$: 

The first subcase is when $s_t'$ does not exist. Then an optimal policy chooses $s_t = \sigma_t$ because in this subcase, a minimizer to the optimization problem defining $J_N(x,t)$ is $u^*_t(x) = 0$ for all $x\in[0,c_{\max}]$. Note that by setting $s_t$ equal to $\sigma_t$, the second policy region vanishes, which ensures that $u^*_t(x) = 0$ for all $x\in[0,c_{\max}]$ with our policy.

The rest of the proof considers the second subcase in which $s_t'$ exists. We set $s_t = s_t'$, since 
\begin{equation}
u^*_t(x) = \begin{cases}
S_t' - x, & \text{if } 0 \leq x < s_t'\\
0, & \text{if } s_t' \leq x\leq c_{\max}
\end{cases}
\end{equation}
is optimal for the optimization problem defining $J_N(x,t)$. We now observe that for $x \in [0,s_t)$, $J_{N}(x,t) = G(S_t',t) + K - c\cdot x$, which is non-increasing in $x$. For $x \in [s_t,c_{\max}]$, $J_{N}(x,t) = G(x,t) - c\cdot x$, which is convex in $x$ since $G(x,t)$ is convex by Proposition 3. Clearly, $J_N(x,t)$ is continuous for $x \in [0,s_t)$ since it is linear in this region, and $J_N(x,t)$ is continuous for $x \in (s_t,c_{\max}]$ since $G(\cdot,t)$ is continuous by Proposition 3. The only question about continuity occurs at $x = s_t$. Since $G(S_t',t) + K = G(s_t,t)$ by definition, we have that $G(S_t',t)+K - c\cdot s_t = G(s_t,t) - c\cdot s_t$. This proves that $J_N(x,t)$ must be continuous at $x = s_t$ since the left and right side of the last equality with $G(\cdot,t)$ are the limits of $J_N(\cdot)$ in the two respective regions. This means $J_N(x,t)$ is continuous on $x\in[0,c_{\max}]$. Because $J_R(\cdot,t)$ is continuous by \textbf{A3}, this means that $J_R(\sigma_t,t) \geq J_N(\sigma_t,t)$ and $J_R(\sigma_t',t) \geq J_N(\sigma_t',t)$ by the definitions of $\sigma_t$ and $\sigma_t'$.

Next, consider any $x\in[\sigma_t,s_t]$ (note that our argument still holds even if this set is empty). Since $J_N(x,t)$ is non-increasing for $x \in [0,s_t]$, and since $J_R(x,t)$ is non-decreasing in $x$ by \textbf{A3}, we have
\begin{equation}
J_R(x,t) \geq J_R(\sigma_t,t) \geq J_N(\sigma_t,t) \geq J_N(x,t)
\end{equation}
for $x\in[\sigma_t,s_t]$. This argument also implies that $J_R(s_t,t) \geq J_N(s_t,t)$. Next, consider any $x\in[s_t,\sigma_t']$, and observe that there exists $\mu\in[0,1]$ such that $x = \mu s_t + (1-\mu)\sigma_t'$. Since we showed above that $J_N(\cdot,t)$ is convex, this means that for any $x\in[s_t,\sigma_t']$, we have 
\begin{equation}
\begin{aligned}
J_N(x,t) = J_N(\mu s_t + (1-\mu)\sigma_t',t) &\leq \mu J_N(s_t,t) + (1-\mu) J_N(\sigma_t',t) \\
&\leq \mu J_R(s_t,t) + (1-\mu) J_R(\sigma_t',t) \\
&\leq J_R(\mu s_t + (1-\mu)\sigma_t',t) = J_R(x,t),
\end{aligned}
\end{equation}
where the last inequality follows because $J_R(\cdot,t)$ is concave by \textbf{A3}. Combining the above shows that $J_R(x,t) \geq J_N(x,t)$ for all $x\in[\sigma_t,\sigma_t']$.

Now, observe that $J_R(x,t) < J_{N}(x,t)$ for $x \in [0,\sigma_t)$ and for $x \in (\sigma_t',c_{\max}]$. (If this last statement were not true, then we could choose an $x' \in [0,\sigma_t)$ or an $x' \in (\sigma_t',c_{\max}]$ such that $J_R(x',t) \geq J_N(x',t)$. Hence, we would have $x' \in I_R$, which reaches a contradiction since, by the definition of $\sigma_t$ and $\sigma_t'$, we would have $\sigma_t \leq x'$ or $\sigma_t' \geq x'$.) Consequently, an optimal policy uses this value of $\sigma_t$. If $\sigma_t' < c_{\max}$, then an optimal policy uses $\Sigma_t = \sigma_t'$, and if $\sigma_t' = c_{\max}$, then an optimal policy chooses $\Sigma_t = +\infty$.

Finally, we must choose a correct value for $S_t$. We must consider three sub-subcases. The first sub-subcase is when $S_t' < \sigma_t$. Since $s_t' \leq S_t'$ by definition of $s_t'$, this means the second policy region is empty when $s_t = s_t'$, which ensures that $u^*_t(x) = 0$ for all $x\in[\sigma_t,\sigma_t']$ with our policy. Thus, we can choose $S_t = S_t'$. (Any arbitrary choice of $S_t$ would give an optimal policy because the corresponding region is empty.) The second sub-subcase is when $\sigma_t \leq S_t' \leq \sigma_t'$. Then an optimal policy chooses $S_t = S_t'$. The third sub-subcase is when $S_t' > \sigma_t'$. By definition of $\sigma_t$ and $\sigma_t'$, this means $J_R(S_t',t) < J_{N}(S_t', t)$ and $\sigma_t' \geq \sigma_t$. Thus, in this sub-subcase, we have
\begin{equation}
J_N(\sigma_t',t) = G(\sigma_t',t) - c\cdot \sigma_t' \geq G(S_t',t) - c\cdot S_t' = J_N(S_t',t),
\end{equation}
where we have used the definition of $S_t'$ in the last inequality. Recalling that $J_R(x,t) \geq J_N(x,t)$ for all $x\in[\sigma_t,\sigma_t']$, we have
\begin{equation}
J_R(S_t',t) < J_N(S_t',t) \leq J_{N}(\sigma_t',t) \leq J_R(\sigma_t',t).
\end{equation}
However, this last statement is a contradiction since $\sigma_t' < S_t'$ and $J_R(\cdot,t)$ is non-decreasing by \textbf{A3}. Therefore, this sub-subcase is not possible. $\Halmos$
\endproof

\begin{figure}[h]
\centering
{\renewcommand{\arraystretch}{0.7}
\begin{tabular}{p{1.4cm} p{1cm} p{2.4cm} p{1cm} p{1cm} p{1cm} p{2.8cm} p{3.6cm}}
\hline
$I_R$ & $s_t'$ & $S_t'$ & $\sigma_t$ & $s_t$ & $S_t$ & $\Sigma_t$ & Interpretation \\
\hline
$I_R = \emptyset$ &  &  & $+\infty$ & $+\infty$ & - & $+\infty$ & Always reset \\
$I_R \neq \emptyset$ & $\nexists s_t'$ &  & $\sigma_t$  & $\sigma_t$ & - & $\sigma_t' \ (\sigma_t' < c_{\max})$, & Missing $[\sigma_t,s_t)$ region \\
 &  &  &  &  &  & $+\infty \ (\sigma_t' = c_{\max})$ &  \\
 & $\exists s_t'$ & $S_t' < \sigma_t$ & $\sigma_t$ & $s_t'$ & - & $\sigma_t'$ or $+\infty$ & Missing $[\sigma_t,s_t)$ region \\
 &  & $\sigma_t \leq S_t' \leq \sigma_t'$ & $\sigma_t$ & $s_t'$ & $S_t'$ & $\sigma_t'$ or $+\infty$ & Four-stage structure \\
 &  & $S_t' > \sigma_t'$ & - & - & - & - & Not possible \\
\hline
\end{tabular}
}
\caption{Proof sketch for Proposition \ref{prop:ps}}
\end{figure}

\proof{Proof of Proposition \ref{prop:vf4s}:} 
The policy structure in Proposition 5 implies that the value function $J(x,t)$ takes the form (\ref{eqn:jtz}). 
This means we have to analyze the derivative of at most four pieces. In the first and fourth regions, by \textbf{A3}, we have that $J(x,t)$ has a derivative that is Lipschitz with constant $L_t$, and that the derivative is absolutely bounded by $\eta_t$. In the second region, $J(x,t)$ is linear in $x$ and thus has a Lipschitz derivative that is absolutely bounded by $c$. In the third region, we know that the first term has a Lipschitz derivative because it is twice differentiable by \textbf{A2}, and that the second term has a Lipschitz derivative by Proposition 2. Thus, the third region has a Lipschitz derivative because it is the sum of two functions with Lipschitz derivatives. To bound its derivative, we note that $\partial_x H(x,t)$ is absolutely bounded by $\kappa_t$ by \textbf{A2}. From Proposition 2, we know that the derivative of the second term is absolutely bounded by $\gamma M_{t+1}$. Thus, the third piece has a derivative that is absolutely bounded by $\kappa_t + \gamma M_{t+1}$. This means the derivative of each piece is absolutely bounded by $\eta_t + c + \kappa_t + \gamma M_{t+1}$. As $0 \leq \gamma \leq m_t/((Lc_{\max}+P)M_{t+1})$, the derivative of each piece is absolutely bounded by $M_t = \eta_t + c + \kappa_t + m_t/(Lc_{\max}+P)$. $\Halmos$
\endproof

\proof{Proof of Theorem \ref{thm:fin}: }
First, note that we can choose $M_k = \eta_k$ using \textbf{A3}, since an optimal policy at $t = k$ is to choose $r_k^*(x) = 1$ for all $x\in[0,c_{\max}]$, which means $J(x,k) = J_0 + \mathbb{E}(s(x,k,w))$. Inductively applying Propositions 5 and 6 implies the desired result. $\Halmos$
\endproof

\proof{Proof of Proposition \ref{prop:sctc}: }
\textbf{A1} and \textbf{A4} hold by the definition of the water problem, and \textbf{A5} holds by assumption. To show \textbf{A2}, we first note that 
\begin{equation}
\begin{aligned}
\xi_n +u_n- (\xi_n+u_n-w_n)^+ &= \xi_n +u_n- (\xi_n + u_n-w_n)^+ - (\xi_n + u_n-w_n)^- + (\xi_n + u_n-w_n)^-\\
&= \xi_n +u_n - (\xi_n + u_n-w_n) + (\xi_n + u_n-w_n)^-\\
&= w_n+ (\xi_n + u_n-w_n)^-
\end{aligned}
\end{equation}
This means the single stage cost can be rewritten as
\begin{equation}
\begin{aligned}
g(\xi_n,\tau_n,u_n,w_n) &= c_u u_n - p\cdot(\xi_n+u_n-w_n)^- + q(\tau_n)\cdot(\xi_n + u_n - (\xi_n+u_n-w_n)^+)\\
&=c_u u_n - (p-q(\tau_n))\cdot(\xi_n+u_n-w_n)^- + q(\tau_n)\cdot w_n
\end{aligned}
\end{equation}
Hence, (\ref{eqn:expsc}) 
is satisfied by setting $K = 0$ and 
\begin{equation}
H(z,t)=-(p-q(t)) \int_z^\infty (z-w) f(w) dw + q(t) \cdot \mathbb{E}(w).
\end{equation}
Note that using the Leibniz integral rule twice gives
\begin{equation}
\label{eqn:derivsh}
\begin{aligned}
\partial_z H(z,t) &= -(p-q(t)) \int_z^\infty f(w) dw \\
\partial_{zz} H(z,t) &= (p-q(t)) \cdot f(z)
\end{aligned}
\end{equation}
Thus, we have $|\partial_z H(z,t)| = |-(p-q(t)) \int_z^\infty f(w) dw| \leq p-q(t)$ since $f(w)$ is a density. Since $H(\cdot,t)$ is twice differentiable, this means it is absolutely continuous. Thus, $H(\cdot,t)$ is strongly convex on $[0,c_{\max}]$ since $\partial_{zz} H(z,t) = (p-q(t)) \cdot f(z) \geq (p-q(t)) \cdot \inf_{z \in [0, c_{\max}]} f(z) = (p-q(t)) \cdot \min_{z \in [0, c_{\max}]} f(z) > 0$ for $z\in[0,c_{\max}]$, where the second equality holds because a Lipschitz continuous function attains its minimum in a compact domain. This shows that \textbf{A2} holds. To show that \textbf{A3} holds, we first note that the water storage problem as described in Section 2.1 implies $\zeta = 0$. Next, note that the reset cost $R(x,t) = c_r x$ is linear. Hence, the conditions of \textbf{A3} follow immediately. $\Halmos$
\endproof

\proof{Proof of Proposition \ref{prop:sctc2}: }
\textbf{A1} and \textbf{A4} hold by the definition of the retail problem, and \textbf{A5} holds by assumption. To show \textbf{A2}, we first note that the single stage cost can be written as
\begin{equation}
\begin{aligned}
g(\xi_n,\tau_n,u_n,w_n) &= c_u u_n + k_u \cdot \mathbf{1}_{\mathbb{R}_{++}}(u_n) - p\cdot(\xi_n+u_n-w_n)^- + q\cdot(\xi_n+u_n-w_n)^+
\end{aligned}
\end{equation}
Hence, (\ref{eqn:expsc}) 
is satisfied by setting $K = k_u$ and 
\begin{equation}
H(z,t)=-p \int_z^\infty (z-w) f(w) dw + q \int_{-\infty}^z (z-w) f(w) dw.
\end{equation}
Note that using the Leibniz integral rule twice gives
\begin{equation}
\begin{aligned}
\partial_z H(z,t) &= -p \int_z^\infty f(w) dw + q \int_{-\infty}^z f(w) dw \\
\partial_{zz} H(z,t) &= (p+q) \cdot f(z)
\end{aligned}
\end{equation}
Thus, we have $|\partial_z H(z,t)| = |-p \int_z^\infty f(w) dw + q \int_{-\infty}^z f(w) dw| = (p+q) \cdot F(z) - p \leq q$, since $F(z) \leq 1$ for $z\in[0,c_{\max}]$. Since $H(\cdot,t)$ is twice differentiable, this means it is absolutely continuous. Thus, $H(\cdot,t)$ is strongly convex on $[0,c_{\max}]$ since $\partial_{zz} H(z,t) = (p+q) \cdot f(z) \geq (p+q) \cdot \min_{z \in [0, c_{\max}]} f(z) > 0$ for $z\in[0,c_{\max}]$. This shows that \textbf{A2} holds. To show \textbf{A3} holds, we first note that the retail management problem as described in Section 2.2 implies $\zeta = 0$. Next, note that the reset cost $R(x,t) = c_r x + k_r$ is linear. Hence, the conditions of \textbf{A3} follow immediately. $\Halmos$
\endproof

\section{Statement and Proof of Lemmas}



\begin{lemma}
\label{lemma:nodec}
Consider two functions $f(x)$ and $g(x)$ defined on a compact interval $[a,b]$. Suppose that $f(x)$ is differentiable, that its derivative is bounded from below by a positive constant $D_x f(x) \geq M > 0$, and that $g(x)$ is Lipschitz continuous with constant $L \geq 0$. Then for any fixed $\lambda$ such that $0 \leq \lambda \leq M/L$, the function $f(x) + \lambda g(x)$ is non-decreasing and continuous.
\end{lemma}

\proof{Proof of Lemma \ref{lemma:nodec}: }
The continuity of $f(x) + \lambda g(x)$ is immediate by the differentiability of $f(x)$ and the Lipschitz continuity of $g(x)$. So we focus on showing $f(x) + \lambda g(x)$ is non-decreasing. Consider any $x_1, x_2 \in[a,b]$ with $x_2 \geq x_1$. Since $D_x f(x) \geq M$, we have $f(x_2) - f(x_1) \geq M(x_2-x_1)$. Since $g(x)$ is Lipschitz continuous with constant $L$, we have $|g(x_2) - g(x_1)| \leq L (x_2 - x_1)$. Combining these two inequalities implies
\begin{equation}
\begin{aligned}
f(x_2) - f(x_1) + \lambda (g(x_2) - g(x_1)) &\geq M(x_2 - x_1) - \lambda L (x_2 - x_1)\\
&\geq (M-\lambda L)\cdot(x_2-x_1) \geq 0.
\end{aligned}
\end{equation}
where the last inequality follows because $x_2 \geq x_1$ and because $0\leq\lambda \leq M/L$, which means $M - \lambda L \geq 0$. Since we have shown $f(x_2) - f(x_1) + \lambda (g(x_2) - g(x_1)) \geq 0$ for any $x_2 \geq x_1$, this implies $f(x) + \lambda g(x)$ is non-decreasing. $\Halmos$
\endproof

\begin{lemma}
\label{lem:conv}
If $\int_{-\infty}^\infty |f(x)|dx$ is finite,  $\int_{-\infty}^\infty |g(x)|dx$ is finite, $f(x) = 0$ for $x < 0$, $|f(x)| \leq M$, $|g(x)| \leq P$, and on $x \geq 0$, we have that $f(x)$ is Lipschitz with constant $L$, then the convolution $(f*g)(x) = \int_{-\infty}^\infty f(x-z)g(z)dz$ is Lipschitz with constant $(L\cdot\textstyle\int_{-\infty}^\infty|g(z)|dz+MP)$.
\end{lemma}

\proof{Proof of Lemma \ref{lem:conv}: }
Without loss of generality, we assume that $x \leq y$. Observe that
\begin{equation}
\begin{aligned}
\textstyle|(f*g)(x)-(f*g)(y)| & = \textstyle|\int_{-\infty}^\infty f(x-z)g(z)dz - \int_{-\infty}^\infty f(y-z)g(z)dz| \\
& = \textstyle|\int_{-\infty}^\infty (f(x-z)-f(y-z))\cdot g(z)dz| \\
& = \textstyle|\int_{-\infty}^x (f(x-z)-f(y-z))\cdot g(z)dz - \int_x^y f(y-z)\cdot g(z) dz| \\
& \leq \textstyle|\int_{-\infty}^x (f(x-z)-f(y-z))\cdot g(z)dz| + |\int_x^y f(y-z)\cdot g(z) dz| \\
& \leq \textstyle\int_{-\infty}^x |f(x-z)-f(y-z)|\cdot |g(z)|dz + \int_x^y|f(y-z)|\cdot |g(z)| dz \\
& \leq \textstyle\int_{-\infty}^x L|(x-z)-(y-z)|\cdot |g(z)|dz + \int_x^y MP\cdot dz \\
& \leq \textstyle\int_{-\infty}^x L|x-y|\cdot |g(z)|dz + MP\cdot|x-y| \\
& \leq (L\cdot\textstyle\int_{-\infty}^\infty |g(z)|dz + MP)\cdot|x-y|\\
\end{aligned}
\end{equation}
This shows the convolution is Lipschitz continuous. $\Halmos$
\endproof

\begin{lemma}
\label{lem:last_lemma}
Consider a function $f(x)$ that is finite at the point $y$ (i.e., $|f(y)|$ is bounded). If $f(x)$ is Lipschitz with constant $L$ on a compact domain $[a,b]$ with $y\in[a,b]$, then it is finitely bounded $|f(x)| \leq L(b-a) + |f(y)|$ for $x\in[a,b]$.
\end{lemma}

\proof{Proof of Lemma \ref{lem:last_lemma}: }
For any $x \in [a,b]$, we have
\begin{equation}
\begin{aligned}
|f(x)| &= |f(x) - f(y) + f(y)| \\
&\leq |f(x) - f(y)| + |f(y)| \\
&\leq L|x-y| + |f(y)| \\
&\leq L(b-a) + |f(y)|
\end{aligned}
\end{equation}
Since the interval $[a,b]$ is compact, this means $a,b$ are finite. This gives the desired bound. $\Halmos$
\endproof

\end{APPENDICES}


\bibliographystyle{informs2014} 
\bibliography{water} 

\begin{thebibliography}{49}
\providecommand{\natexlab}[1]{#1}
\providecommand{\url}[1]{\texttt{#1}}
\providecommand{\urlprefix}{URL }

\bibitem[{Arrow et~al.(1951)Arrow, Harris, \protect\BIBand{}
  Marschak}]{arrow1951}
Arrow KJ, Harris T, Marschak J (1951) Optimal inventory policy.
  \emph{Econometrica} 19(3):250--272.

\bibitem[{Barnes \protect\BIBand{} Lea-Greenwood(2006)}]{barnes2006}
Barnes L, Lea-Greenwood G (2006) Fast fashioning the supply chain: shaping the
  research agenda. \emph{Journal of Fashion Marketing and Management}
  10(3):259--271.

\bibitem[{Bertsekas(1995)}]{bertsekas1995}
Bertsekas DP (1995) \emph{Dynamic programming and optimal control}, volume~1
  (Athena Scientific).

\bibitem[{Bhardwaj \protect\BIBand{} Fairhurst(2010)}]{bhardwaj2010}
Bhardwaj V, Fairhurst A (2010) Fast fashion: response to changes in the fashion
  industry. \emph{The International Review of Retail, Distribution and Consumer
  Research} 20(1):165--173.

\bibitem[{Blackburn \protect\BIBand{} Scudder(2009)}]{blackburn2009}
Blackburn J, Scudder G (2009) Supply chain strategies for perishable products:
  the case of fresh produce. \emph{Production and Operations Management}
  18(2):129--37.

\bibitem[{Botev(2017)}]{botev2017}
Botev ZI (2017) The normal law under linear restrictions: simulation and
  estimation via minimax tilting. \emph{Journal of the Royal Statistical
  Society: Series B (Statistical Methodology)} 79(1):125--148.

\bibitem[{Cachon \protect\BIBand{} Swinney(2011)}]{cachon2011}
Cachon GP, Swinney R (2011) The value of fast fashion: Quick response, enhanced
  design, and strategic consumer behavior. \emph{Management Science}
  57(4):778--795.

\bibitem[{Caro \protect\BIBand{} Gallien(2007)}]{caro2007}
Caro F, Gallien J (2007) Dynamic assortment with demand learning for seasonal
  consumer goods. \emph{Management Science} 53(2):276--292.

\bibitem[{Caro \protect\BIBand{} Gallien(2010)}]{caro2010}
Caro F, Gallien J (2010) Inventory management of a fast-fashion retail network.
  \emph{Operations Research} 58(2):257--273.

\bibitem[{Caro \protect\BIBand{}
  {Mart{\'i}nez-de-Alb{\'e}niz}(2015)}]{caro2015}
Caro F, {Mart{\'i}nez-de-Alb{\'e}niz} V (2015) Fast fashion: Business model
  overview and research opportunities. Agrawal N, Smith SA, eds., \emph{Retail
  Supply Chain Management}, 237--264, International Series in Operations
  Research \& Management Science (Springer).

\bibitem[{Chu et~al.(2005)Chu, Hsu, \protect\BIBand{} Shen}]{chu2005}
Chu LY, Hsu VN, Shen Z (2005) An economic lot‐sizing problem with perishable
  inventory and economies of scale costs: Approximation solutions and worst
  case analysis. \emph{Naval Research Logistics} 52(6):536--548.

\bibitem[{Coelho \protect\BIBand{} Laporte(2014)}]{coelho2014}
Coelho LC, Laporte G (2014) Optimal joint replenishment, delivery and inventory
  management policies for perishable products. \emph{Computers \& Operations
  Research} 47:42--52.

\bibitem[{Coelho et~al.(2003)Coelho, James, Sunna, Abu~Jaish, \protect\BIBand{}
  Chatila}]{coelho2003}
Coelho ST, James S, Sunna N, Abu~Jaish A, Chatila J (2003) Controlling water
  quality in intermittent supply systems. \emph{Water Supply} 3(1-2):119--125.

\bibitem[{de~Farias \protect\BIBand{} Van~Roy(2003)}]{de2003linear}
de~Farias DP, Van~Roy B (2003) The linear programming approach to approximate
  dynamic programming. \emph{Operations Research} 51(6):850--865.

\bibitem[{Farahani et~al.(2012)Farahani, Grunow, \protect\BIBand{}
  G{\"u}nther}]{farahani2012}
Farahani P, Grunow M, G{\"u}nther HO (2012) Integrated production and
  distribution planning for perishable food products. \emph{Flexible Services
  and Manufacturing Journal} 24:28--51.

\bibitem[{Fremlin(2001)}]{fremlin2001}
Fremlin DH (2001) \emph{Measure Theory}, volume~2 (Torres Fremlin).

\bibitem[{Gadgil(1998)}]{gadgil1998}
Gadgil A (1998) Drinking water in developing countries. \emph{Annual Review of
  Energy and the Environment} 23(1):253--286.

\bibitem[{Hadley \protect\BIBand{} Whitin(1963)}]{hadley1963}
Hadley G, Whitin T (1963) \emph{Analysis of inventory systems} (Prentice-Hall,
  Inc.).

\bibitem[{Haskell et~al.(2016)Haskell, Jain, \protect\BIBand{}
  Kalathil}]{haskell2016empirical}
Haskell WB, Jain R, Kalathil D (2016) Empirical dynamic programming.
  \emph{Mathematics of Operations Research} 41(2):402--429.

\bibitem[{Hsu(2000)}]{hsu2000}
Hsu VN (2000) Dynamic economic lot size model with perishable inventory.
  \emph{Management Science} 46(8):1159--1169.

\bibitem[{Hsu \protect\BIBand{} Lowe(2001)}]{hsu2001}
Hsu VN, Lowe TJ (2001) Dynamic economic lot size models with
  period-pair-dependent backorder and inventory costs. \emph{Operations
  Research} 49(2):316--321.

\bibitem[{Karaesmen et~al.(2011)Karaesmen, {Scheller-Wolf}, \protect\BIBand{}
  Deniz}]{karaesmen2011}
Karaesmen IZ, {Scheller-Wolf} A, Deniz B (2011) Managing perishable and aging
  inventories: Review and future research directions. Kempf KG, Keskinocak P,
  Uzsoy R, eds., \emph{Planning Production and Inventories in the Extended
  Enterprise: A State of the Art Handbook}, 393--436 (Springer).

\bibitem[{Kariotoglou et~al.(2013)Kariotoglou, Summers, Summers, Kamgarpour,
  \protect\BIBand{} Lygeros}]{kariotoglou2013approximate}
Kariotoglou N, Summers S, Summers T, Kamgarpour M, Lygeros J (2013) Approximate
  dynamic programming for stochastic reachability. \emph{2013 European Control
  Conference (ECC)}, 584--589.

\bibitem[{Klingel(2012)}]{klingel2012}
Klingel P (2012) Technical causes and impacts of intermittent water
  distribution. \emph{Water Supply} 12(4):504--512.

\bibitem[{Kumpel \protect\BIBand{} Nelson(2016)}]{kumpel2016}
Kumpel E, Nelson KL (2016) Intermittent water supply: Prevalence, practice, and
  microbial water quality. \emph{Environmental Science \& Technology}
  50(2):542--553.

\bibitem[{Lee \protect\BIBand{} Schwab(2005)}]{lee2005}
Lee EJ, Schwab KJ (2005) Deficiencies in drinking water distribution systems in
  developing countries. \emph{Journal of Water and Health} 3(2):109--127.

\bibitem[{Maier et~al.(2009)Maier, Pepper, \protect\BIBand{} Gerba}]{maier2009}
Maier RM, Pepper IL, Gerba CP (2009) \emph{Environmental Microbiology}
  (Elsevier).

\bibitem[{Matsuo \protect\BIBand{} Ogawa(2007)}]{matsuo2007}
Matsuo H, Ogawa S (2007) Innovating innovation: The case of seven-eleven japan.
  \emph{International Commerce Review} 7(2).

\bibitem[{Mintz et~al.(2017)Mintz, Shen, \protect\BIBand{} Aswani}]{mintz2017}
Mintz Y, Shen Z, Aswani A (2017) Local water storage control for the developing
  world. \emph{2017 IEEE 56th Annual Conference on Decision and Control (CDC)},
  5074--5079.

\bibitem[{Nagurney et~al.(2013)Nagurney, Yu, Masoumi, \protect\BIBand{}
  Nagurney}]{nagurney2013}
Nagurney A, Yu M, Masoumi AH, Nagurney LS (2013) \emph{Networks against time:
  Supply chain analytics for perishable products} (Springer).

\bibitem[{Nahmias(1982)}]{nahmias1982}
Nahmias S (1982) Perishable inventory theory: A review. \emph{Operations
  Research} 30(4):680--708.

\bibitem[{Pierskalla(2005)}]{pierskalla2005}
Pierskalla WP (2005) Supply chain management of blood banks. \emph{Operations
  research and health care}, 103--145 (Springer).

\bibitem[{Pierskalla \protect\BIBand{} Roach(1972)}]{pierskalla1972}
Pierskalla WP, Roach CD (1972) Optimal issuing policies for perishable
  inventory. \emph{Management Science} 18(11):603--614.

\bibitem[{Powell(2007)}]{powell2007approximate}
Powell WB (2007) \emph{Approximate Dynamic Programming: Solving the curses of
  dimensionality} (John Wiley \& Sons).

\bibitem[{Prastacos(1984)}]{prastacos1984}
Prastacos GP (1984) Blood inventory management: an overview of theory and
  practice. \emph{Management Science} 30(7):777--800.

\bibitem[{Rosenfield(1989)}]{rosenfield1989}
Rosenfield DB (1989) Disposal of excess inventory. \emph{Operations Research}
  37(3):404--409.

\bibitem[{Rosenfield(1992)}]{rosenfield1992}
Rosenfield DB (1992) Optimality of myopic policies in disposing excess
  inventory. \emph{Operations Research} 40(4):800--803.

\bibitem[{Ryzhov et~al.(2012)Ryzhov, Powell, \protect\BIBand{}
  Frazier}]{ryzhov2012knowledge}
Ryzhov IO, Powell WB, Frazier PI (2012) The knowledge gradient algorithm for a
  general class of online learning problems. \emph{Operations Research}
  60(1):180--195.

\bibitem[{Scarf(1959)}]{scarf1959}
Scarf H (1959) The optimality of {(S, s)} policies in the dynamic inventory
  problem. Technical report, Stanford University.

\bibitem[{Shen et~al.(2011)Shen, Dessouky, \protect\BIBand{}
  Ord{\'o}{\~n}ez}]{shen2011}
Shen Z, Dessouky M, Ord{\'o}{\~n}ez F (2011) Perishable inventory management
  system with a minimum volume constraint. \emph{Journal of the Operational
  Research Society} 62(12):2063--2082.

\bibitem[{Snyder \protect\BIBand{} Shen(2019)}]{snyder2019}
Snyder LV, Shen Z (2019) \emph{Fundamentals of Supply Chain Theory} (Wiley).

\bibitem[{Tokajian \protect\BIBand{} Hashwa(2003)}]{tokajian2003}
Tokajian S, Hashwa F (2003) Water quality problems associated with intermittent
  water supply. \emph{Water Science \& Technology} 47(3):229--234.

\bibitem[{{United Nations}(2015)}]{ohchr2015}
{United Nations} (2015) The right to water. Technical report, Office of the
  United Nations High Commissioner for Human Rights.

\bibitem[{{van Donselaar} et~al.(2006){van Donselaar}, {van Woensel},
  Broekmeulen, \protect\BIBand{} Fransoo}]{van2006}
{van Donselaar} K, {van Woensel} T, Broekmeulen R, Fransoo J (2006) Inventory
  control of perishables in supermarkets. \emph{International Journal of
  Production Economics} 104(2):462--472.

\bibitem[{WaterAid(2016)}]{wateraid2016}
WaterAid (2016) Water: At what cost? the state of the world's water 2016.
  Technical report, WaterAid.

\bibitem[{WHO \protect\BIBand{} UNICEF(2017)}]{who2017}
WHO, UNICEF (2017) Progress on drinking water, sanitation and hygiene: 2017
  update and {SDG} baselines. Technical report, WHO/UNICEF Joint Monitoring
  Programme for Water Supply, Sanitation and Hygiene.

\bibitem[{{World Bank}(2010)}]{worldbank2010}
{World Bank} (2010) The {Karnataka} urban water sector improvement project: 24
  x 7 water supply is achievable. Technical report, Water and Sanitation
  Program.

\bibitem[{Z{\'e}rah(1998)}]{zerah1998}
Z{\'e}rah MH (1998) How to assess the quality dimension of urban
  infrastructure: The case of water supply in {Delhi}. \emph{Cities}
  15(4):285--290.

\bibitem[{Zidek \protect\BIBand{} Kolmanovsky(2016)}]{zidek2016stochastic}
Zidek RA, Kolmanovsky IV (2016) Stochastic drift counteraction optimal control
  and enhancing convergence of value iteration. \emph{2016 IEEE 55th Conference
  on Decision and Control (CDC)}, 1119--1124.

\end{thebibliography}


\end{document}